\documentclass[11pt,reqno]{amsart}
\usepackage[T1]{fontenc}
\usepackage{lmodern}
\usepackage{amsmath,amssymb,mathtools}
\usepackage[margin=1.05in]{geometry}
\usepackage{microtype}
\usepackage{enumitem}
\usepackage[colorlinks=true,linkcolor=blue,citecolor=blue,urlcolor=blue]{hyperref}
\hypersetup{pdftitle={Stability of Hilbert--Samuel functions with linear bounds in homological degree},pdfauthor={Ton That Quoc Tan and Van Duc Trung},pdfsubject={Commutative algebra},pdfkeywords={Hilbert function, perturbation, homological degree, arithmetic degree, deficiency module, equicharacteristic local ring}}
\numberwithin{equation}{section}
\newtheorem{theorem}{Theorem}[section]
\newtheorem{lemma}[theorem]{Lemma}
\newtheorem{proposition}[theorem]{Proposition}
\newtheorem{corollary}[theorem]{Corollary}
\theoremstyle{definition}
\newtheorem{example}[theorem]{Example}
\newcommand{\m}{\mathfrak m}
\newcommand{\q}{\mathfrak q}
\newcommand{\f}{\mathbf f}
\newcommand{\hdeg}{\operatorname{hdeg}}
\newcommand{\adeg}{\operatorname{adeg}}
\newcommand{\Ann}{\operatorname{Ann}}
\newcommand{\Ass}{\operatorname{Ass}}
\newcommand{\Ext}{\operatorname{Ext}}
\setlist[enumerate,1]{label=\textup{(\roman*)},leftmargin=2.2em}
\title[Stability of Hilbert--Samuel functions]{Stability of Hilbert--Samuel functions\\with linear bounds in homological degree}
\author{Ton That Quoc Tan}
\address{Department of Mathematics, FPT University, Da Nang, Vietnam}
\email{tanttq@fe.edu.vn, quoctanmath@gmail.com}
\author{Van Duc Trung}
\address{Department of Mathematics, University of Education, Hue University, 34 Le Loi, Hue, Vietnam}
\email{vdtrung@hueuni.edu.vn}
\date{}
\begin{document}
\begin{abstract}
Let $(R,\m)$ be a Noetherian local ring and let $I$ be generated by a filter regular sequence. Set $A=R/I$, $s=\dim A>0$, and $D=\hdeg(A)$. If $A$ is equicharacteristic, we prove that perturbations of order $(2s+1)D$ preserve the Hilbert--Samuel function in every dimension. In arbitrary characteristic, we obtain the order $[s+(s+1)(k!)^{k-1}]D^k$, where $k$ is the maximum of $1$ and the dimensions of the lower deficiency modules of $\widehat A$. Both bounds depend only on the final quotient. The proof constructs monic annihilating equations for the deficiency modules and uses them to choose a parameter $d$-sequence with an explicit relation to a reduction of the maximal ideal. In the equicharacteristic case, a dimension filtration and characteristic polynomials over regular subrings give equations of degree at most the arithmetic degree, with coefficients in the required powers of the reduction.
\end{abstract}
\maketitle

\section{Introduction}

All rings are commutative and Noetherian, and all modules are finitely
generated unless stated otherwise. We write $\ell_R(-)$ for $R$-module
length and omit the subscript when the base ring is clear.
The ideal generated by an empty sequence is understood to be zero.
Let $(R,\m)$ be a local ring and let $\f=(f_1,\ldots,f_r)$ be a sequence
in $\m$, with $r\ge1$. Put $I=(\f)=(f_1,\ldots,f_r)R$ and $A=R/I$.
For an integer $N\ge1$, a small perturbation of $\f$ by a tuple
$\boldsymbol\varepsilon=(\varepsilon_1,\ldots,\varepsilon_r)$ replaces $I$ by
\[
 I_{\boldsymbol\varepsilon}=(f_1+\varepsilon_1,\ldots,f_r+\varepsilon_r),
 \qquad \varepsilon_i\in\m^N.
\]
For any proper ideal $J$ of $R$, write
\[
 H_{R/J}(n)=\ell_R\bigl(R/(J+\m^n)\bigr)\qquad(n\ge1).
\]
In particular, $H_A(n)=\ell_R(A/\m^nA)$.
The problem is to find an explicit order $N$ for which $H_{R/I_{\boldsymbol\varepsilon}}=H_A$ for every such perturbation. This quantitative form of finite determinacy also describes how accurately defining equations must be approximated to preserve the Hilbert--Samuel function.
To state the consequences and compare examples, we use the Hilbert
perturbation index
\[
 p_{\m}(\f)=\min\{N\ge1:H_{R/I_{\boldsymbol\varepsilon}}=H_A
       \text{ for all }\varepsilon_i\in\m^N\},
\]
with value $\infty$ when no such integer exists.

Srinivas and Trivedi initiated this study in
\cite[Theorem~3 and Corollary~5]{ST96} and obtained an explicit
multiplicity bound in the Cohen--Macaulay case in
\cite[Proposition~1]{ST97}. Ma, Quy and Smirnov proved stability for
filter regular sequences over arbitrary local rings
\cite[Theorem~3.5]{MQS}. Quy and Ngo Viet Trung obtained explicit
bounds involving Artin--Rees numbers and colon modules
\cite[Theorem~3.5]{QT}. The bounds in \cite[Theorem~4.5]{QT} involve degrees of the
ambient ring and intermediate quotients. For generalized
Cohen--Macaulay ambient rings, final-quotient bounds are given in
\cite[Corollary~4.8]{QT}, while Quy and Van Duc Trung
\cite[Theorem~1.4]{QVT} give a Hilbert bound involving the final
homological degree and the Buchsbaum invariant of the ambient ring.

We seek a bound determined by the final quotient alone.
The homological degree is taken with respect to the maximal ideal and
is defined recursively in Section~\ref{sec:module-annihilators}.
A hat denotes completion with respect to the maximal ideal. For a finite
module $M$ over a local ring $R$ that need not be complete, our convention is
\[
 \hdeg_R(M):=\hdeg_{\widehat R}(M\otimes_R\widehat R).
\]
We write $e(A)$ for the Hilbert--Samuel multiplicity of $A$ with respect
to its maximal ideal. Recall that $f_1,\ldots,f_r$ is filter regular if
each module
\[
 \frac{(f_1,\ldots,f_{i-1}):_R f_i}{(f_1,\ldots,f_{i-1})}
 \qquad(1\le i\le r)
\]
has finite length. A local ring is equicharacteristic if its
characteristic equals that of its residue field, or equivalently if
it contains a field. Our principal result gives a linear bound in every
positive dimension for equicharacteristic final quotients.

\begin{theorem}[Equicharacteristic final quotients]
\label{thm:equicharacteristic}
Let $(R,\m)$ be a Noetherian local ring and let $\f=(f_1,\ldots,f_r)$
be a filter regular sequence in $\m$, with $r\ge1$.
Put $I=(\f)$ and $A=R/I$. Suppose that $A$ is equicharacteristic and
$s=\dim A>0$, and put $D=\hdeg(A)$. Then
\begin{equation}\label{eq:equicharacteristic-main}
 p_{\m}(\f)\le (2s+1)D.
\end{equation}
The same bound holds if filter regularity is replaced by
\begin{equation}\label{eq:koszul-hypothesis}
 \ell_R(H_j(\f;R))<\infty\qquad\text{for every }j>0.
\end{equation}
\end{theorem}

Theorem~\ref{thm:equicharacteristic} gives an affirmative answer to
the case of Quy and Ngo Viet Trung's \cite[Question~4.9]{QT} in which
the filtration is maximal-ideal-adic, the extended degree is the
homological degree, and the final quotient is equicharacteristic.

The equicharacteristic hypothesis covers localizations of finitely
generated algebras over fields and quotients of formal power series
rings over fields, in both characteristic zero and positive
characteristic. Neither $R$ nor $A$ is assumed to be generalized
Cohen--Macaulay, and $A$ need not be reduced or equidimensional.
The hypothesis concerns only $A$; the ambient ring $R$ may have
mixed characteristic.

Without a restriction on characteristic, we obtain a polynomial bound
whose exponent is controlled by the lower deficiency modules.
For a finite module $M$ over a complete local ring $(B,\mathfrak n)$,
define
\[
 K^i(M):=\operatorname{Hom}_B
 \bigl(H^i_{\mathfrak n}(M),E_B(B/\mathfrak n)\bigr)
 \qquad(i\in\mathbb Z),
\]
where $H^i$ denotes local cohomology and $E_B(B/\mathfrak n)$ is
the injective hull of the residue field over $B$.
For $A=R/I$, we apply this definition to
$\widehat A=\widehat R/I\widehat R$, regarded as a finite
$\widehat R$-module. Set
\begin{equation}\label{eq:k-definition}
 k(A)=\max\bigl(\{1\}\cup
       \{\dim K^i(\widehat A):0\le i<s,\ K^i(\widehat A)\ne0\}\bigr).
\end{equation}
Since $\dim K^i(\widehat A)\le i$, we have
\[
 1\le k(A)\le\max\{1,s-1\}.
\]

\begin{theorem}[Arbitrary characteristic]\label{thm:main}
Let $(R,\m)$ be a Noetherian local ring and let $\f=(f_1,\ldots,f_r)$
be a filter regular sequence in $\m$, with $r\ge1$.
Put $I=(\f)$ and $A=R/I$, and suppose that $s=\dim A>0$.
Put $D=\hdeg(A)$, $k=k(A)$, and
\[
 C_{s,k}=s+(s+1)(k!)^{k-1}.
\]
Then
\[
 p_{\m}(\f)\le C_{s,k}D^k.
\]
The same bound holds if filter regularity is replaced by
\eqref{eq:koszul-hypothesis}.
\end{theorem}

In both theorems, $H_j(\f;R)$ denotes Koszul homology. Condition
\eqref{eq:koszul-hypothesis} is equivalent to regularity of the
localized sequence at every nonmaximal prime containing $I$.
Indeed, at such a prime the Koszul criterion identifies vanishing
of positive homology with regularity; at a prime not containing $I$,
one element is a unit and the complex is contractible. A finite
module has finite length precisely when its support is contained
in the singleton consisting of the maximal ideal. A filter regular
sequence satisfies this condition because its successive colon
modules vanish at every nonmaximal prime. We use the Koszul
hypothesis to replace the generators one at a time and to test
filter regularity on the quotient by the other generators.
The conclusion concerns the Hilbert function of the final quotient;
it does not require preservation of filter regularity of the
perturbed ordered sequence.

Both bounds depend only on $A$. The general bound is linear when
all lower deficiency modules have dimension at most one, which
always holds for $s\leq2$. Section~\ref{sec:consequences} records
further linear consequences and compares our bounds with the
homological-degree bound in \cite[Theorem~4.5]{QT}. In Example~\ref{ex:qt-hdeg-comparison},
the final quotient is fixed, has dimension three and satisfies
$k(A)=2$, so it lies outside Corollary~\ref{cor:linear}, but
Theorem~\ref{thm:equicharacteristic} applies. Our sufficient order
is constant, while the cited bound grows with the homological
degree of the ambient ring.
The linear bound also holds in mixed characteristic when every lower
deficiency module is annihilated by a power of the residue characteristic;
see Corollary~\ref{cor:prime-power-deficiencies}.

After passing to a complete local ring with infinite residue field,
we construct monic annihilating equations for the action of elements
of the maximal ideal on a finite module $M$ of positive dimension $d$.
For an equation of degree $a$, the coefficient of $Z^{a-j}$
lies in the $j$-th power of an ideal generated by $d$ general
elements, for $1\le j\le a$. In equal characteristic, the degree
is $\adeg(M)\le\hdeg(M)$: we adapt Vasconcelos's filtration and
characteristic-polynomial argument
\cite[Proposition~6 and Theorem~7]{Vasconcelos1996}
to complete local rings, using the order valuation to obtain the
required coefficient bounds. In arbitrary characteristic,
homological-degree estimates, the Eakin--Sathaye theorem, and a
determinant argument give degree at most $(d!)^{d-1}\hdeg(M)^d$.

We apply these equations to the deficiency modules of $A$.
Starting from a reduction $Q\subseteq R$ of $\m$ on $A$, we
construct an ideal $\q\subseteq R$ whose image in $A$ is generated
by a triangular system of parameters. Annihilator propagation through
local cohomology makes this system a $d$-sequence, so
$\operatorname{reg}\operatorname{gr}_{\q}(A)=0$.
The triangular form also yields an exact identity between powers
of $Q$ and $\q$ in $R$.

Finally, we combine stability for one filter regular element with
successive generic replacements of the generators allowed by the
Koszul hypothesis. This preserves the $\q$-adic Hilbert function,
while regularity zero and Huneke's intersection theorem allow
comparison with the maximal-ideal filtration. The resulting
perturbation orders depend only on the final quotient.

Section~\ref{sec:module-annihilators} develops the module estimates.
Section~\ref{sec:equicharacteristic-annihilators} proves the
linear-degree monic bound in equal characteristic.
Section~\ref{sec:parameters} constructs the parameter ideals.
Section~\ref{sec:transport} proves the Hilbert-function comparison under
perturbations. Section~\ref{sec:main-proof} derives a common
perturbation estimate and gives separate proofs of the two main theorems. Section~\ref{sec:consequences} gives the consequences
and a comparison example.

\section{Homological degree and monic annihilators}
\label{sec:module-annihilators}

Throughout this section, $(B,\m)$ is a complete Noetherian local ring
with infinite residue field $\Bbbk=B/\m$, and $M$ denotes a finite
$B$-module. We write
\[
 \mu(M)=\mu_B(M)=\dim_{\Bbbk}(M/\m M)
\]
for its minimal number of generators, and $\ell(M)=\ell_B(M)$ for its
length. Multiplicities and homological degrees are taken with respect
to $\m$. For $\dim M=d>0$, the
multiplicity $e(M)$ is $d!$ times the leading coefficient of the
polynomial agreeing with $\ell(M/\m^nM)$ for large $n$.
An ideal $Q\subseteq\m$ is a reduction of $\m$ on $M$ if
$\m^{a+1}M=Q\m^aM$ for some $a\geq0$.

We use the deficiency modules $K^i(M)$ defined in the Introduction.
By local duality, these modules are finite, $\dim K^i(M)\leq i$, and
$K^i(M)=0$ for $i<0$ or $i>\dim M$. If a complete regular
local ring $S$ of dimension $v$ maps onto $B$, then
\[
 K^i(M)\cong\operatorname{Ext}^{v-i}_S(M,S).
\]
The Ext description also explains the dimension bound directly.
If $\mathfrak p\in\operatorname{Spec}S$ has
$\operatorname{ht}\mathfrak p<v-i$, then the regular local ring
$S_{\mathfrak p}$ has global dimension $\operatorname{ht}\mathfrak p$,
and hence
$\operatorname{Ext}^{v-i}_{S_{\mathfrak p}}
 (M_{\mathfrak p},S_{\mathfrak p})=0$.
Thus every prime in the support of $K^i(M)$ has height at least
$v-i$, which gives $\dim K^i(M)\le i$.
In particular, $K^0(M)$ has length $\ell(H^0_{\m}(M))$ by
Matlis duality. The homological degree is defined over any complete
local ring by the following recursion.

Throughout the paper, we use the conventions $\dim 0=-\infty$ and $\hdeg(0)=0$.  For a
nonzero module of dimension $d$, the homological degree is defined
recursively by
\begin{equation}
 \hdeg(M)=
 \begin{cases}
  \ell(M),&d=0,\\[2mm]
  e(M)+\displaystyle\sum_{i=0}^{d-1}
          \binom{d-1}{i}\hdeg(K^i(M)),&d>0.
 \end{cases}
 \label{eq:module-hdeg-definition}
\end{equation}
The recursion is well defined because $\dim K^i(M)<d$ in the sum.
These invariants may equally be computed over a quotient of $B$
that acts on $M$. Indeed, let $C=B/J$, where $J\subseteq\m$ and $JM=0$. Local cohomology
with respect to $\m$ agrees with that with respect to $\m C$, and
$E_C(\Bbbk)\cong\operatorname{Hom}_B(C,E_B(\Bbbk))$.
Hence adjunction gives
\[
 \operatorname{Hom}_B(H^i_{\m}(M),E_B(\Bbbk))
 \cong
 \operatorname{Hom}_C(H^i_{\m C}(M),E_C(\Bbbk)).
\]
Thus the deficiency modules agree, as do their dimensions over $B$
and $C$. Lengths and maximal-ideal filtrations also agree, so
multiplicities and, by recursion, homological degrees are the same.
For a Cohen--Macaulay module, $\hdeg(M)=e(M)$.  The invariant was
introduced and developed in \cite[Definition~2.8 and Section~3]{Vasconcelos1998}.

For an ideal $\mathfrak a$ of $B$, we write
\[
 \operatorname{gr}_{\mathfrak a}(M)
 =\bigoplus_{n\ge0}\mathfrak a^nM/\mathfrak a^{n+1}M,
\]
viewed as a graded module over $\operatorname{gr}_{\mathfrak a}(B)$.
For a nonnegatively graded ring $G=\bigoplus_{n\ge0}G_n$, we write
$G_+=\bigoplus_{n>0}G_n$ for its positive-degree ideal.
Throughout the paper, if a module is defined over a quotient of the
ambient ring, an ideal appearing as a subscript in local cohomology or
an associated graded object refers to its image in that quotient.
We also use the same symbols for elements and their images when the
ambient ring is clear.

The word \emph{general} refers to a nonempty Zariski open subset of the
appropriate finite-dimensional vector space over $\Bbbk$.  Thus a
general element of $\m$ means any lift of a class in a suitable open
subset of $\m/\m^2$.
Since $\Bbbk$ is infinite, every finite-dimensional $\Bbbk$-vector
space is irreducible in the Zariski topology. Hence any finite
intersection of nonempty Zariski open subsets of the same vector
space is nonempty. We use this observation throughout to impose
finitely many general conditions simultaneously. In particular, a
general element is superficial and filter regular on any prescribed
finite collection of positive-dimensional finite modules.  Here
filter regular means that its annihilator on the module has finite
length.

We use superficiality in its usual filtration sense: $x\in\m$
is superficial on $M$ if, for some $c$,
$(\m^{n+1}M:_M x)\cap\m^cM=\m^nM$ for all $n\geq c$.
Such general choices exist because the degree-one initial form
can avoid the finitely many associated primes of
$\operatorname{gr}_{\m}(M)$ that do not contain its positive-degree
ideal. Each avoidance condition excludes a proper linear subspace.
Filter regularity on $M$ is obtained by avoiding
$\Ass(M)\setminus\{\m\}$; these also give proper subspaces of
$\m/\m^2$, since $\m\subseteq\mathfrak p+\m^2$ would imply
$\m=\mathfrak p$ by Nakayama's lemma in $B/\mathfrak p$.

We will also use the following observation over any local ring.
If $x$ is filter regular on $M$ and $T=H^0_{\m}(M)$, then $x$
is regular on $M/T$. Indeed, if $xm\in T$, some power of $\m$
sends $m$ into $0:_Mx\subseteq T$, and a further power kills it;
hence $m\in T$. Cancellation on $M/T$ then gives
\begin{equation}\label{eq:filter-regular-torsion}
 T\cap x^qM=x^qT\qquad(q\geq1).
\end{equation}

The next lemma is the local form of
\cite[Proposition~2.14]{Vasconcelos1998}, stated there for graded
modules. We include the local argument.

\begin{lemma}
\label{lem:hdeg-finite-length}
Let $0\longrightarrow L\longrightarrow M\longrightarrow N
\longrightarrow0$ be an exact sequence of finite $B$-modules.
\begin{enumerate}
 \item If $\ell(L)<\infty$, then
       $\hdeg(M)=\ell(L)+\hdeg(N)$.
 \item If $\ell(N)<\infty$, then
       $\hdeg(M)\leq\hdeg(L)+\ell(N)$.
\end{enumerate}
\end{lemma}

\begin{proof}
For the given exact sequence,
$\operatorname{Supp}M=\operatorname{Supp}L\cup\operatorname{Supp}N$,
so $\dim M=\max\{\dim L,\dim N\}$. If $M$ has finite length, so
do $L$ and $N$, and both assertions follow from additivity of length.
We may therefore assume $d=\dim M>0$.

For \textup{(i)}, the finite length of $L$ gives $\dim N=d$.
Since $H^i_{\m}(L)=0$ for $i>0$, the local cohomology sequence gives
\[
 0\longrightarrow L\longrightarrow H^0_{\m}(M)
 \longrightarrow H^0_{\m}(N)\longrightarrow0,
 \qquad H^i_{\m}(M)\cong H^i_{\m}(N)\quad(i\ge1).
\]
Matlis duality is exact and preserves the lengths of finite-length
modules. Thus
\[
 \ell(K^0(M))=\ell(L)+\ell(K^0(N)),\qquad
 K^i(M)\cong K^i(N)\quad(i\ge1).
\]
Also $e(M)=e(N)$ by additivity of multiplicity, because
$\dim L<d$. Since $K^0(M)$ and $K^0(N)$ have finite length,
substitution into \eqref{eq:module-hdeg-definition} proves
\textup{(i)}.

For \textup{(ii)}, the finite length of $N$ gives $\dim L=d$ and
$e(L)=e(M)$, again by additivity of multiplicity. The local
cohomology sequence is
\[
 \begin{aligned}
 0\longrightarrow H^0_{\m}(L)\longrightarrow H^0_{\m}(M)
 \longrightarrow N\longrightarrow H^1_{\m}(L)
 \longrightarrow H^1_{\m}(M)\longrightarrow0,
 \end{aligned}
\]
and $H^i_{\m}(L)\cong H^i_{\m}(M)$ for $i\ge2$.
Write $(-)^{\vee}=\operatorname{Hom}_B(-,E_B(\Bbbk))$.
The injectivity of $E_B(\Bbbk)$ makes this contravariant functor
exact. Applying it to the displayed sequence yields
\[
 0\longrightarrow K^1(M)\longrightarrow K^1(L)
 \longrightarrow N^{\vee}\longrightarrow K^0(M)
 \longrightarrow K^0(L)\longrightarrow0,
\]
while $K^i(M)\cong K^i(L)$ for $i\ge2$.
Since $\ell(N^{\vee})=\ell(N)$, exactness gives
\[
 \ell(K^0(M))\leq\ell(K^0(L))+\ell(N).
\]
If $d=1$, this inequality and $e(M)=e(L)$ already prove
\textup{(ii)}, because the defining formula for the homological
degree then contains only the $K^0$ term.

Suppose $d\ge2$. The displayed dual sequence identifies $K^1(M)$
with a submodule of $K^1(L)$ whose cokernel embeds into
$N^{\vee}$ and hence has finite length. Both modules have dimension
at most one. If $K^1(L)$ has finite length, including the case
$K^1(L)=0$, then
$\hdeg(K^1(M))=\ell(K^1(M))\le\ell(K^1(L))=\hdeg(K^1(L))$.
Otherwise $\dim K^1(L)=1$. Its finite-length quotient by $K^1(M)$
then implies $\dim K^1(M)=1$ and
$e(K^1(M))=e(K^1(L))$. Moreover, the inclusion of the modules
restricts to
\[
 H^0_{\m}(K^1(M))\subseteq H^0_{\m}(K^1(L)).
\]
For every finite module $U$ of dimension one,
\eqref{eq:module-hdeg-definition} and Matlis duality give
$\hdeg(U)=e(U)+\ell(H^0_{\m}(U))$.
The equality of multiplicities and the last inclusion therefore give
\[
 \hdeg(K^1(M))\leq\hdeg(K^1(L))
\]
in this case as well. Finally, the defining formula, the equality
$e(M)=e(L)$, and the isomorphisms for $i\ge2$ give
\[
 \hdeg(M)\le e(L)+\sum_{i=0}^{d-1}\binom{d-1}{i}\hdeg(K^i(L))+\ell(N)
           =\hdeg(L)+\ell(N).
\]
This proves \textup{(ii)}.
\end{proof}

The following is the local form of
\cite[Corollary~2.16]{Vasconcelos1998}; see also its use for local
rings in the proof of \cite[Theorem~3.5]{Vasconcelos1998}.

\begin{lemma}
\label{lem:generic-power-hdeg}
Let $M$ be a finite $B$-module and fix an integer $q\ge1$.
For a general element $x\in\m$, we have
\begin{equation}
 \hdeg(M/x^qM)\leq q\hdeg(M).
 \label{eq:generic-power-hdeg}
\end{equation}
\end{lemma}

\begin{proof}
Fix an integer $q\ge1$. We prove the assertion by induction on
$d=\dim M$.
The assertion is immediate for $M=0$. If $d=0$, then
\[
 \hdeg(M/x^qM)=\ell(M/x^qM)\le\ell(M)\le q\hdeg(M).
\]

Suppose $d>0$. Put $T=H^0_{\m}(M)$ and $U=M/T$.
Since $\ell(T)<\infty$, we have $\dim U=d$; the local cohomology
sequence of $0\to T\to M\to U\to0$ gives $H^0_{\m}(U)=0$.
Choose $x\in\m$ general for $M$, $U$, and the modules $K^i(U)$,
$0\le i<d$. We may take $x$ superficial and filter regular on
these modules and, by the induction hypothesis, require
\[
 \hdeg(K^i(U)/x^qK^i(U))\le q\hdeg(K^i(U))
 \qquad(0\le i<d).
\]
These conditions can be imposed simultaneously, since only finitely many
modules are involved and $\dim K^i(U)\le i<d$. In particular, $x$ is regular
on $U$. By \eqref{eq:filter-regular-torsion},
\[
 0\longrightarrow T/x^qT\longrightarrow M/x^qM
 \longrightarrow U/x^qU\longrightarrow0
\]
is exact. Lemma~\ref{lem:hdeg-finite-length}\textup{(i)} gives
\[
 \hdeg(M/x^qM)=\ell(T/x^qT)+\hdeg(U/x^qU),
 \qquad \hdeg(M)=\ell(T)+\hdeg(U).
\]
Thus it suffices to prove $\hdeg(U/x^qU)\le q\hdeg(U)$,
since $\ell(T/x^qT)\le\ell(T)\le q\ell(T)$.

Set $V=U/x^qU$. We first compute its dimension and multiplicity.
Choose $c$ such that
$(\m^{n+1}U:_U x)\cap\m^cU=\m^nU$ for $n\ge c$.
Artin--Rees applied to $xU\subseteq U$ gives
$\m^{n+1}U\cap xU\subseteq x\m^cU$ for $n\gg0$.
Since $x$ is regular on $U$, this implies
$(\m^{n+1}U:_U x)\subseteq\m^cU$ and hence
$(\m^{n+1}U:_U x)=\m^nU$ for $n\gg0$.
Iterating this equality $q$ times gives
$(\m^nU:_U x^q)=\m^{n-q}U$ for $n\gg0$.
Therefore
\[
 0\longrightarrow U/\m^{n-q}U
 \xrightarrow{x^q}U/\m^nU
 \longrightarrow V/\m^nV\longrightarrow0
\]
is exact for $n\gg0$. If $P_U(n)=\ell(U/\m^nU)$ for $n\gg0$,
then the Hilbert polynomial of $V$ is $P_U(n)-P_U(n-q)$.
It has degree $d-1$ and leading coefficient $qe(U)/(d-1)!$.
Consequently, $\dim V=d-1$ and
\[
 \ell(V)=qe(U)\quad\text{if }d=1,
 \qquad e(V)=qe(U)\quad\text{if }d\ge2.
\]
If $d=1$, then $U$ is Cohen--Macaulay, so
$\hdeg(V)=\ell(V)=qe(U)=q\hdeg(U)$.

Suppose $d\ge2$. The local cohomology sequence of
\[
 0\longrightarrow U\xrightarrow{x^q}U
 \longrightarrow V\longrightarrow0
\]
gives, for every $j\ge0$,
\[
 0\longrightarrow H^j_{\m}(U)/x^qH^j_{\m}(U)
 \longrightarrow H^j_{\m}(V)
 \longrightarrow 0:_{H^{j+1}_{\m}(U)}x^q
 \longrightarrow0.
\]
Matlis duality is exact and interchanges the kernel and cokernel
of multiplication by $x^q$. Hence
\begin{equation}
 0\longrightarrow K^{j+1}(U)/x^qK^{j+1}(U)
 \longrightarrow K^j(V)
 \longrightarrow 0:_{K^j(U)}x^q\longrightarrow0.
 \label{eq:power-deficiency-sequence}
\end{equation}
For $0\le j\le d-2$, filter regularity makes $x$ regular on
$K^j(U)/H^0_{\m}(K^j(U))$. Thus
\[
 0:_{K^j(U)}x^q\subseteq H^0_{\m}(K^j(U)),\qquad
 \ell(0:_{K^j(U)}x^q)
 \le\ell(H^0_{\m}(K^j(U)))\le\hdeg(K^j(U)).
\]
Write $h_i=\hdeg(K^i(U))$.
Lemma~\ref{lem:hdeg-finite-length}\textup{(ii)}, applied to
\eqref{eq:power-deficiency-sequence}, and the induction hypothesis
give
\[
 \hdeg(K^j(V))
 \le\hdeg(K^{j+1}(U)/x^qK^{j+1}(U))
       +\ell(0:_{K^j(U)}x^q)
 \le qh_{j+1}+h_j\le q(h_{j+1}+h_j).
\]
Using $e(V)=qe(U)$ and Pascal's identity, we obtain
\[
 \begin{aligned}
 \hdeg(V)
 &\le qe(U)+q\sum_{j=0}^{d-2}\binom{d-2}{j}(h_j+h_{j+1})\\
 &=q\left(e(U)+\sum_{i=0}^{d-1}\binom{d-1}{i}h_i\right)
 =q\hdeg(U).
 \end{aligned}
\]
Together with the preceding reduction from $M$ to $U$, this proves
\eqref{eq:generic-power-hdeg}.
\end{proof}

\begin{proposition}[Sally-type estimate for submodules]
\label{prop:module-sally-bound}
Every finite $B$-module $N$ satisfies $\mu(N)\leq\hdeg(N)$.
Let $M$ have dimension $d\geq1$, and let $q\geq1$ be an integer. If a
submodule $L\subseteq M$ contains $\m^qM$, then
\begin{equation}
 \mu(L)\leq q^{d-1}\hdeg(M)+(d-1)\mu(M).
 \label{eq:module-sally-bound}
\end{equation}
In particular,
\begin{equation}
 \mu(\m^qM)\leq\bigl(q^{d-1}+d-1\bigr)\hdeg(M).
 \label{eq:module-power-generators}
\end{equation}
\end{proposition}

\begin{proof}
For the first assertion, the case of finite length follows from
$\mu(N)\leq\ell(N)=\hdeg(N)$. For $\dim N>0$, we argue by
induction on $\dim N$. Choose a general
$x\in\m$ satisfying Lemma~\ref{lem:generic-power-hdeg} with $q=1$.
We may also take $x$ filter regular on $N$, so it avoids the minimal
primes of $\operatorname{Supp}N$ of dimension $\dim N$. Thus
$\dim N/xN=\dim N-1$, and induction gives
\[
 \mu(N)=\mu(N/xN)\leq\hdeg(N/xN)\leq\hdeg(N).
\]
The equality holds because $xN\subseteq\m N$.
See also \cite[Corollary~3.3]{Vasconcelos1998}.

We now prove the submodule estimate by induction on $d=\dim M\geq1$, following the
argument of \cite[Theorem~3.5]{Vasconcelos1998} for ideals.
First let $d=1$. Set $T=H^0_{\m}(M)$, $U=M/T$, and let
$\overline L$ denote the image of $L$ in $U$.
Since $\dim U=1$ and $H^0_{\m}(U)=0$, the module $U$ is
Cohen--Macaulay. Moreover, $\m^qU\subseteq\overline L$, so
$U/\overline L$ has finite length. Consequently, $\overline L$ has
dimension one and $e(\overline L)=e(U)=e(M)$.
A nonzerodivisor on $U$
is also one on its submodule $\overline L$; hence $\overline L$
is Cohen--Macaulay. The inequality just proved gives
\[
 \mu(\overline L)\leq\hdeg(\overline L)=e(M).
\]
The sequence $0\to L\cap T\to L\to\overline L\to0$ now gives
\[
 \mu(L)\leq\mu(L\cap T)+\mu(\overline L)
 \leq\ell(T)+e(M)=\hdeg(M).
\]

For $d\geq2$, choose a general $x\in\m$ that is filter regular on
$M$ and satisfies Lemma~\ref{lem:generic-power-hdeg} for the given
$q$. Since $x^qM\subseteq\m^qM\subseteq L$, the modules
$\overline M=M/x^qM$ and $\overline L=L/x^qM$ satisfy
$\overline L\subseteq\overline M$ and
$\m^q\overline M\subseteq\overline L$. Moreover,
$\operatorname{Supp}\overline M=\operatorname{Supp}(M/xM)$, so
$\dim\overline M=d-1$.
The exact sequence
\[
 0\longrightarrow x^qM\longrightarrow L\longrightarrow
 \overline L\longrightarrow0
\]
gives the first inequality below. Induction on $d$, followed by
Lemma~\ref{lem:generic-power-hdeg}, then gives
\begin{align*}
 \mu(L)
 &\leq\mu(\overline L)+\mu(x^qM)\\
 &\leq q^{d-2}\hdeg(\overline M)
          +(d-2)\mu(\overline M)+\mu(M)\\
 &\leq q^{d-1}\hdeg(M)+(d-1)\mu(M).
\end{align*}
Here $\mu(\overline M)=\mu(M)$ since $x^qM\subseteq\m M$, and
$\mu(x^qM)\leq\mu(M)$ since $x^qM$ is an image of $M$.
Taking $L=\m^qM$ and using $\mu(M)\leq\hdeg(M)$ proves
the last assertion.
\end{proof}

The following is a module version of
\cite[Corollary~3.6]{Vasconcelos1998}.
We adapt Vasconcelos's argument, combining
Proposition~\ref{prop:module-sally-bound} with the
Eakin--Sathaye theorem \cite[Theorem~1]{EakinSathaye1976}.
We give the details of the passage to modules.

\begin{corollary}[Reduction bound for modules]
\label{cor:module-reduction-bound}
Let $M$ have dimension $d\geq1$, and put $H=\hdeg(M)$ and
$n=d!H$.  For general elements $y_1,\ldots,y_d\in\m$, the ideal
$Q=(y_1,\ldots,y_d)$ satisfies
\begin{equation}
 \m^nM=Q\m^{n-1}M.
 \label{eq:module-reduction-equality}
\end{equation}
\end{corollary}

\begin{proof}
By Proposition~\ref{prop:module-sally-bound},
\[
 \mu(\m^nM)\leq H(n^{d-1}+d-1)<\binom{n+d}{d}.
\]
For the strict inequality, if $d=1$, then $n=H$ and the assertion
is $H<n+1$. If $d\geq2$, expansion of the product defining the
binomial coefficient gives
\begin{align*}
 \binom{n+d}{d}
 &>\frac{n^d}{d!}+\frac{d(d+1)}{2d!}n^{d-1}\\
 &=Hn^{d-1}+\frac{d(d+1)}2Hn^{d-2}\\
 &\geq H(n^{d-1}+d-1).
\end{align*}
The first inequality is strict because the omitted terms are positive.
The equality uses $n=d!H$, and the last inequality follows from
$n^{d-2}\geq1$ and $d(d+1)/2\geq d-1$.

Put $G=\operatorname{gr}_{\m}(B)$ and
$F=\operatorname{gr}_{\m}(M)$. The $G$-module $F$ is generated
by $F_0=M/\m M$, and
$\dim_{\Bbbk}F_n=\mu(\m^nM)$.
Choose nonzero elements $v_1,\ldots,v_r\in F_0$ generating $F$.
Then
\[
 F_t=\sum_{j=1}^r(Gv_j)_t\qquad(t\geq0).
\]
For each $j\in\{1,\ldots,r\}$, multiplication by $v_j$ gives a
graded isomorphism
\[
 A_j:=G/\Ann_G(v_j)\ \cong\ Gv_j.
\]
Since $v_j\ne0$ and $G_0=\Bbbk$ is a field, $(A_j)_0=\Bbbk$;
thus $A_j$ is a standard graded $\Bbbk$-algebra.
The inclusion $Gv_j\subseteq F$ gives
\[
 \dim_{\Bbbk}(A_j)_n
 \leq\dim_{\Bbbk}F_n<\binom{n+d}{d}.
\]

Fix $j$. We first find suitable linear forms in $A_j$.
Set $C=(A_j)_{(A_j)_+}$ and $\mathfrak b=(A_j)_+C$.
Since $A_j$ is standard graded,
\[
 (A_j)_+^t=\bigoplus_{q\geq t}(A_j)_q,\qquad
 (A_j)_+^t/(A_j)_+^{t+1}\cong(A_j)_t
 \qquad(t\geq0).
\]
Each quotient is annihilated by $(A_j)_+$, so localization at
$(A_j)_+$ does not change it. Consequently
\[
 \mathfrak b^t/\mathfrak b^{t+1}\cong(A_j)_t,\qquad
 \mu_C(\mathfrak b^n)=\dim_{\Bbbk}(A_j)_n<\binom{n+d}{d}.
\]
By \cite[Theorem~1]{EakinSathaye1976}, there are
$h_{1,j},\ldots,h_{d,j}\in\mathfrak b$ such that
$\mathfrak b^n=\sum_{i=1}^d h_{i,j}\mathfrak b^{n-1}$.
Let $\lambda_{i,j}\in(A_j)_1$ correspond to
$h_{i,j}+\mathfrak b^2$ under
$\mathfrak b/\mathfrak b^2\cong(A_j)_1$; this class may be zero.
Reducing
$\mathfrak b^n=\sum_{i=1}^d h_{i,j}\mathfrak b^{n-1}$
modulo $\mathfrak b^{n+1}$ and using
$\mathfrak b^t/\mathfrak b^{t+1}\cong(A_j)_t$ gives
\[
 (A_j)_n=\sum_{i=1}^d\lambda_{i,j}(A_j)_{n-1}.
\]

The forms $\lambda_{i,j}$ may depend on $j$. We now choose a
single tuple in $G_1^d$ that works for every $j$.
For each $j$, let $\pi_j:G_1\to(A_j)_1$ be the quotient map.
For each tuple $\boldsymbol\ell=(\ell_1,\ldots,\ell_d)\in G_1^d$,
define
\[
 \Phi_{j,\boldsymbol\ell}:(A_j)_{n-1}^{\oplus d}\longrightarrow(A_j)_n,
 \qquad
 (u_1,\ldots,u_d)\longmapsto
 \sum_{i=1}^d\pi_j(\ell_i)u_i.
\]
Since $\pi_j$ is surjective, choose
$\ell_i^{(j)}\in G_1$ with
$\pi_j(\ell_i^{(j)})=\lambda_{i,j}$ for $1\leq i\leq d$.
For $\boldsymbol\ell^{(j)}
=(\ell_1^{(j)},\ldots,\ell_d^{(j)})$, we have
\[
 \operatorname{Im}\Phi_{j,\boldsymbol\ell^{(j)}}
 =\sum_{i=1}^d\lambda_{i,j}(A_j)_{n-1}
 =(A_j)_n.
\]
Thus, for each $j$, the map $\Phi_{j,\boldsymbol\ell}$
is surjective for at least one tuple $\boldsymbol\ell$.

Fix a basis of $G_1$, and choose bases of $(A_j)_{n-1}$ and
$(A_j)_n$ for each $j$.
The entries of the matrix of $\Phi_{j,\boldsymbol\ell}$ are linear
polynomials in the coordinates of $\boldsymbol\ell$.
If $(A_j)_n\ne0$, surjectivity at one tuple implies that some
minor of size $\dim_{\Bbbk}(A_j)_n$ is a nonzero polynomial;
denote one such minor by $\Delta_j$.
If $(A_j)_n=0$, put $\Delta_j=1$, since the map is then
surjective for every tuple.
The product
\[
 \Delta=\prod_{j=1}^r\Delta_j
\]
is a nonzero polynomial in the coordinates of $G_1^d$.
Since $\Bbbk$ is infinite, there are tuples with
$\Delta(\boldsymbol\ell)\ne0$. This condition therefore holds
for general $\ell_1,\ldots,\ell_d\in G_1$ and makes every
$\Phi_{j,\boldsymbol\ell}$ surjective.

Under $A_j\cong Gv_j$, these simultaneous surjections give
\[
 (Gv_j)_n=\sum_{i=1}^d\ell_i(Gv_j)_{n-1}
 \qquad(1\leq j\leq r).
\]
Hence the same forms satisfy
\[
 F_n
 =\sum_{j=1}^r(Gv_j)_n
 =\sum_{j=1}^r\sum_{i=1}^d\ell_i(Gv_j)_{n-1}
 =\sum_{i=1}^d\ell_i\left(\sum_{j=1}^r(Gv_j)_{n-1}\right)
 =\sum_{i=1}^d\ell_iF_{n-1}.
\]
Lift $\ell_i\in G_1=\m/\m^2$ to $y_i\in\m$ and put
$Q=(y_1,\ldots,y_d)$. By the definition of multiplication on $F$,
the last equality says
\[
 \frac{\m^nM}{\m^{n+1}M}
 =\frac{Q\m^{n-1}M+\m^{n+1}M}{\m^{n+1}M},
 \qquad\text{so}\qquad
 \m^nM=Q\m^{n-1}M+\m^{n+1}M.
\]
For the finite module $T=\m^nM/Q\m^{n-1}M$, this is $T=\m T$.
Nakayama's lemma gives $T=0$, proving
\eqref{eq:module-reduction-equality}.
Multiplication by powers of $\m$ gives the equality in every
later degree.
\end{proof}

\begin{proposition}[Monic annihilators in arbitrary characteristic]
\label{prop:monic-module-annihilator}
Let $M$ be a finite $B$-module. If $d=\dim M\geq1$, set
\[
 H=\hdeg(M),\qquad c_d=(d!)^{d-1}.
\]
For general $y_1,\ldots,y_d\in\m$, put $Q=(y_1,\ldots,y_d)$.
There exists an integer $a$ with
\[
 1\leq a\leq c_dH^d,
\]
such that every $z\in\m$ admits a monic polynomial $P_z(Z)\in B[Z]$
\begin{equation}
 P_z(Z)=Z^a+\alpha_{1,z}Z^{a-1}+\cdots+\alpha_{a,z},
 \qquad \alpha_{j,z}\in Q^j,
 \qquad P_z(z)M=0.
 \label{eq:monic-module-annihilator}
\end{equation}
If $M$ has dimension zero, the coefficient and annihilation conditions
hold with $a=\ell(M)$ and $P_z(Z)=Z^a$, for any ideal $Q$.
For $M=0$, one may use the constant polynomial $1$ of degree zero.
\end{proposition}

\begin{proof}
If $M$ has finite length $h$, the strict decrease of the nonzero terms
of its $\m$-adic filtration, by Nakayama's lemma, implies $\m^hM=0$.
Thus $P_z(Z)=Z^h$ has the required coefficient and annihilation
properties for every $z\in\m$ and any ideal $Q$.
This also covers $M=0$, when $h=0$ and $P_z(Z)=1$.

Now assume $d=\dim M\geq1$ and put $H=\hdeg(M)$. Set
\[
 T=H^0_{\m}(M),\qquad h=\ell(T),\qquad U=M/T,
 \qquad H'=\hdeg(U)=H-h.
\]
The last equality follows from
Lemma~\ref{lem:hdeg-finite-length}\textup{(i)}. The module $U$ has
dimension $d$, satisfies $H^0_{\m}(U)=0$, and has $H'\geq1$.
Set $n=d!H'$ and $E=\m^{n-1}U$.
By Corollary~\ref{cor:module-reduction-bound}, general
$y_1,\ldots,y_d\in\m$ give $Q=(y_1,\ldots,y_d)$ with
\[
 \m E=QE.
\]
The module $E$ is nonzero: otherwise a power of $\m$ would
annihilate $U$, contradicting $\dim U=d>0$. Thus $b=\mu(E)\geq1$.

If $n=1$, then $d=H'=1$ and $E=U$, so the first assertion of
Proposition~\ref{prop:module-sally-bound} gives
$b=\mu(U)\leq H'=1=c_d(H')^d$.
If $n\geq2$, Proposition~\ref{prop:module-sally-bound}, applied
with $q=n-1$, gives
\[
 b\leq H'\bigl((n-1)^{d-1}+d-1\bigr)
 \leq H'n^{d-1}=c_d(H')^d.
\]
The second inequality follows from $(n-1)^p+p\leq n^p$ for
integers $n\geq2$ and $p\geq0$. Thus $1\leq b\leq c_d(H')^d$
in all cases.

Choose $b$ generators $u_1,\ldots,u_b$ of $E$. For each
$z\in\m$, the inclusion $zE\subseteq QE$ allows us to write
\[
 zu_i=\sum_{j=1}^b a_{ij}u_j,\qquad a_{ij}\in Q.
\]
Taking the determinant of $ZI_b-(a_{ij})$ gives the monic polynomial
\[
 F_z(Z)=\det(ZI_b-(a_{ij}))
       =Z^b+\beta_{1,z}Z^{b-1}+\cdots+\beta_{b,z},
\]
where $I_b$ is the $b\times b$ identity matrix.
The coefficient of $Z^{b-j}$ is a signed sum of products of $j$
entries of $(a_{ij})$, so $\beta_{j,z}\in Q^j$ for $1\leq j\leq b$.
Applying the adjugate of $zI_b-(a_{ij})$ to these relations gives
$F_z(z)E=0$. It follows that
\[
 \m^{n-1}F_z(z)U=0,
\]
so $F_z(z)U\subseteq H^0_{\m}(U)=0$. Consequently,
$F_z(z)M\subseteq T$.

Since $\m^hT=0$, the polynomial
\[
 P_z(Z)=Z^hF_z(Z)
\]
annihilates $M$ after evaluation at $z$. Put $a=h+b$.
The coefficient of $Z^{a-j}$ in $P_z(Z)$ is
$\beta_{j,z}\in Q^j$ for $1\leq j\leq b$ and is zero for
$b<j\leq a$. Thus all the required coefficient conditions hold.
The fixed degree $a$ satisfies
\[
 1\leq a=h+b\leq h+c_d(H-h)^d\leq c_dH^d.
\]
For the last inequality, use $c_d\geq1$ and
$H^d-(H-h)^d\geq h$ for integers $H\geq1$ and $0\leq h<H$.
\end{proof}

\section{Arithmetic degree and monic annihilators in equal characteristic}
\label{sec:equicharacteristic-annihilators}

We retain the notation and the completeness and residue-field
assumptions of Section~\ref{sec:module-annihilators}. We first recall
the arithmetic degree and a coefficient estimate over a regular local
ring. The equicharacteristic hypothesis enters in
Proposition~\ref{prop:equicharacteristic-annihilator} through the
existence of regular power series subrings.

The underlying degrees in each dimension go back to
Hartshorne~\cite[Chapter~2]{Hartshorne1966}; the terminology follows
Bayer and Mumford~\cite[Section~3]{BayerMumford1993}.

For a nonzero finite $B$-module $M$ of dimension $d$, let
$D_i$ be its largest submodule of dimension at most $i$, and set
$D_{-1}=0$. The resulting dimension filtration is
\[
 0=D_{-1}\subseteq D_0\subseteq\cdots\subseteq D_d=M.
\]
Set $C_i=D_i/D_{i-1}$ and put $e_i=e(C_i)$ if $C_i\ne0$ and
$e_i=0$ otherwise, interpreting multiplicity as length in dimension
zero. We define the arithmetic degree of $M$ by
\begin{equation}\label{eq:arithmetic-degree}
 \adeg(M):=\sum_{i=0}^d e_i
 =\sum_{\mathfrak p\in\Ass_B M}
   \ell_{B_{\mathfrak p}}
     \bigl(H^0_{\mathfrak pB_{\mathfrak p}}(M_{\mathfrak p})\bigr)
     e(B/\mathfrak p).
\end{equation}
We set $\adeg(0)=0$.

To verify the second equality, recall the associativity formula
for a nonzero finite $B$-module $N$:
\[
 e(N)=
 \sum_{\substack{\mathfrak p\in\Ass_B N\\
                  \dim(B/\mathfrak p)=\dim N}}
 \ell_{B_{\mathfrak p}}(N_{\mathfrak p})\,e(B/\mathfrak p).
\]
The dimension filtration satisfies
\[
 \Ass_B C_i
 =\{\mathfrak p\in\Ass_B M:\dim(B/\mathfrak p)=i\},
\]
so each nonzero $C_i$ has dimension $i$.
At every prime $\mathfrak p$ in this set, we have
$(D_{i-1})_{\mathfrak p}=0$, since $\dim D_{i-1}<i$.
Moreover, $(D_i)_{\mathfrak p}$ is the largest finite-length
submodule of $M_{\mathfrak p}$. Hence
\[
 (C_i)_{\mathfrak p}\cong(D_i)_{\mathfrak p}
 =H^0_{\mathfrak pB_{\mathfrak p}}(M_{\mathfrak p}).
\]
Applying the associativity formula to each nonzero $C_i$ and
summing over $i=0,\ldots,d$ proves the second equality in
\eqref{eq:arithmetic-degree}, since these sets partition $\Ass_B M$.

\begin{lemma}\label{lem:arithmetic-homological-degree}
For every finite $B$-module $M$, one has
$\adeg(M)\le\hdeg(M)$.
\end{lemma}
\begin{proof}
The assertion is immediate for a finite-length module, including
the zero module. Suppose $d=\dim M>0$ and use the notation above.
The equality $e_d=e(M)$ follows from $\dim D_{d-1}<d$.
For $0\le i<d$, local duality and the associativity formula give
\[
 e_i=
 \begin{cases}
 e(K^i(M)),&\dim K^i(M)=i,\\
 0,&\dim K^i(M)<i,
 \end{cases}
\]
where multiplicity in dimension zero means length.
To verify this identity, choose a complete regular presentation
$P\twoheadrightarrow B$, with $v=\dim P$.
For $\mathfrak p\in\operatorname{Spec}B$ with
$\dim B/\mathfrak p=i$, let $\mathfrak P$ be its inverse image
in $P$. Then $\dim P_{\mathfrak P}=v-i$, and local duality over
$P_{\mathfrak P}$ identifies
\[
 K^i(M)_{\mathfrak p}
 \cong\Ext^{v-i}_{P_{\mathfrak P}}
             (M_{\mathfrak p},P_{\mathfrak P})
\]
with the Matlis dual of
$H^0_{\mathfrak pB_{\mathfrak p}}(M_{\mathfrak p})$.
For every $\mathfrak p$ with $\dim(B/\mathfrak p)=i$, the
isomorphism $(C_i)_{\mathfrak p}\cong
H^0_{\mathfrak pB_{\mathfrak p}}(M_{\mathfrak p})$ and local duality give
\[
 \ell_{B_{\mathfrak p}}\bigl(K^i(M)_{\mathfrak p}\bigr)
 =\ell_{B_{\mathfrak p}}\bigl((C_i)_{\mathfrak p}\bigr)<\infty.
\]
A nonzero finite-length $B_{\mathfrak p}$-module has
$\mathfrak pB_{\mathfrak p}$ as its only associated prime.
Thus, by localization of associated primes, for every such
$\mathfrak p$ we have
\[
 \mathfrak p\in\Ass_B K^i(M)
 \Longleftrightarrow
 \ell_{B_{\mathfrak p}}\bigl(K^i(M)_{\mathfrak p}\bigr)>0
 \Longleftrightarrow
 \ell_{B_{\mathfrak p}}\bigl((C_i)_{\mathfrak p}\bigr)>0
 \Longleftrightarrow\mathfrak p\in\Ass_B C_i.
\]
Since $\dim(B/\mathfrak p)=i$ for every
$\mathfrak p\in\Ass_B C_i$, it follows that
\[
 \{\mathfrak p\in\Ass_B K^i(M):\dim(B/\mathfrak p)=i\}
 =\Ass_B C_i.
\]
Recall that $\dim K^i(M)\le i$.
If $\dim K^i(M)=i$, the associativity formula gives
$e(K^i(M))=e_i$ by summing these equal local lengths times
$e(B/\mathfrak p)$ over the common set of primes. If
$\dim K^i(M)<i$, this set is empty, so $C_i=0$ and $e_i=0$.
This is the local arithmetic-degree computation underlying
\cite[Proposition~5]{Vasconcelos1996}.
In particular, $e_i\le\hdeg(K^i(M))$. Since
$\binom{d-1}{i}\ge1$, the defining formula for $\hdeg$ yields
\[
 \adeg(M)=e(M)+\sum_{i=0}^{d-1}e_i
 \le e(M)+\sum_{i=0}^{d-1}\binom{d-1}{i}\hdeg(K^i(M))
 =\hdeg(M).
\]
\end{proof}

\begin{lemma}
\label{lem:characteristic-coefficients}
Let $(S,\mathfrak n)$ be a regular local ring of positive dimension,
with fraction field $K=\operatorname{Frac}(S)$. Let $E$ be a finite
torsion-free $S$-module
of rank $e>0$, and let $\varphi$ be an $S$-linear endomorphism of $E$.
Suppose that a submodule $L\subseteq E$ has rank $e$ and satisfies
$\varphi(L)\subseteq\mathfrak nL$. Then the characteristic polynomial
of $\varphi_K:=\varphi\otimes_S\operatorname{id}_K$ on
$E\otimes_SK$ has the form
\[
 F(Z)=Z^e+\beta_1Z^{e-1}+\cdots+\beta_e,
 \qquad \beta_j\in\mathfrak n^j,
\]
and $F(\varphi)E=0$.
\end{lemma}
\begin{proof}
For every height-one prime $\mathfrak p$ of $S$, the ring
$S_{\mathfrak p}$ is regular local of dimension one, hence a
principal ideal domain; see
\cite[Theorems~11.2 and~19.3]{Matsumura1989}.
Since $E_{\mathfrak p}$ is finite and torsion-free over
$S_{\mathfrak p}$, it is flat and hence free of rank $e$ by
\cite[Theorems~7.7 and~7.10]{Matsumura1989}.
Let $A_{\mathfrak p}\in M_e(S_{\mathfrak p})$ represent the localized
endomorphism $\varphi_{\mathfrak p}$ in an $S_{\mathfrak p}$-basis of
$E_{\mathfrak p}$. Under the natural identification
$E_{\mathfrak p}\otimes_{S_{\mathfrak p}}K\cong E\otimes_SK$,
this basis induces a $K$-basis of $E\otimes_SK$, with respect to
which $\varphi_K$ is represented by the same matrix $A_{\mathfrak p}$.
Consequently,
\[
 F(Z)=\det(ZI_e-A_{\mathfrak p})\in S_{\mathfrak p}[Z].
\]
The ring $S$ is normal by \cite[Theorem~19.4]{Matsumura1989}, and
\cite[Proposition~4.10.4]{HunekeSwanson2006} gives
$S=\bigcap_{\operatorname{ht}\mathfrak p=1}S_{\mathfrak p}$ inside $K$.
Since each coefficient of $F$ belongs to $S_{\mathfrak p}$ for every
height-one prime $\mathfrak p$ of $S$, the intersection formula gives
$F\in S[Z]$.
Cayley--Hamilton on $E\otimes_SK$, together with the injection
$E\hookrightarrow E\otimes_SK$, gives $F(\varphi)E=0$.

Since $\operatorname{gr}_{\mathfrak n}(S)$ is a polynomial ring
over $S/\mathfrak n$, the order function
$v(a)=\max\{q:a\in\mathfrak n^q\}$ for $0\ne a\in S$
extends by $v(a/b)=v(a)-v(b)$ to a discrete valuation on $K$;
see \cite[Theorem~6.7.9]{HunekeSwanson2006}.
Set $v(0)=+\infty$. The corresponding valuation ring is
\[
 V=\{u\in K:v(u)\ge0\}.
\]
Then $S\subseteq V\subseteq K$.
Choose $\pi\in\mathfrak n\setminus\mathfrak n^2$, so that $v(\pi)=1$.
For every $a\in\mathfrak n$, we have $v(a/\pi)\ge0$, hence
$a\in\pi V$. Together with $\pi\in\mathfrak n$, this gives
$\mathfrak nV=\pi V$. For $a\in S$ and $j\ge0$, the condition
$a\in\pi^jV$ is equivalent to $v(a)\ge j$, and hence to
$a\in\mathfrak n^j$. Thus
\[
 S\cap\pi^jV=\mathfrak n^j\qquad(j\ge0).
\]
Since $E$ is torsion-free, we identify each $\ell\in L\subseteq E$
with $\ell\otimes 1\in E\otimes_S K$. Let $VL$ denote the
$V$-submodule of $E\otimes_S K$ generated by $L$. Thus
\[
 VL=\left\{\sum_{i=1}^{r} v_i(\ell_i\otimes 1)
      \;\middle|\; r\ge1,\ \ell_i\in L,\ v_i\in V\right\},
\]
where $v_i(\ell_i\otimes 1)=\ell_i\otimes v_i$.
The module $VL$ is finite and torsion-free over
the discrete valuation ring $V$, hence free. Fix a $V$-basis $w_1,\ldots,w_t$ of $VL$.
The vectors $w_1,\ldots,w_t$ are linearly independent over
$K=\operatorname{Frac}(V)$, since clearing denominators in a
$K$-linear relation among the $w_i$ gives a $V$-linear relation.
Since $L\subseteq VL=\sum_{i=1}^t Vw_i$, the subspace
$\operatorname{span}_K\{w_1,\ldots,w_t\}$ contains $L$.
As $L\subseteq E$ and both have rank $e$ over $S$, $L$ spans
$E\otimes_S K$ over $K$. Thus $w_1,\ldots,w_t$ form a $K$-basis
of $E\otimes_S K$ and $t=e$. The hypothesis gives $\varphi_K(VL)\subseteq\pi VL$,
so the matrix $A_V$ of $\varphi_K$ in this basis has entries in $\pi V$.
Hence the coefficient $\beta_j$ of $Z^{e-j}$ in
$F(Z)=\det(ZI_e-A_V)$ belongs to $\pi^jV$. Since $\beta_j\in S$,
the displayed intersection gives $\beta_j\in\mathfrak n^j$.
\end{proof}

\begin{proposition}[Monic annihilators of arithmetic degree]
\label{prop:equicharacteristic-annihilator}
Assume that $B$ is equicharacteristic. Let $M$ be a nonzero finite
$B$-module of dimension $d\ge1$. For general
$y_1,\ldots,y_d\in\m$, put $Q=(y_1,\ldots,y_d)$ and
$a=\adeg(M)$. For every $z\in\m$, there is a monic
polynomial
\[
 P_z(Z)=Z^a+\alpha_{1,z}Z^{a-1}+\cdots+\alpha_{a,z}\in B[Z],
 \qquad \alpha_{j,z}\in Q^j,
 \qquad P_z(z)M=0.
\]
For a nonzero finite-length $B$-module $M$, take
$P_z(Z)=Z^{\ell(M)}$ and any ideal $Q\subseteq B$.
For $M=0$, take $P_z(Z)=1$.
\end{proposition}
\begin{proof}
The finite-length assertion follows as in
Proposition~\ref{prop:monic-module-annihilator}.
Suppose $d\ge1$, and let $C_i$ and $e_i$ be the factors and
multiplicities of the dimension filtration above.
By Corollary~\ref{cor:module-reduction-bound}, general
$y_1,\ldots,y_d$ may be chosen so that, for every nonzero $C_i$
with $i>0$, the ideal $Q_i=(y_1,\ldots,y_i)B$ is a reduction
of $\m$ on $C_i$. Indeed, these conditions define finitely many
nonempty Zariski open subsets of the common parameter space of
$d$-tuples. Their intersection is therefore nonempty, so a general
tuple satisfies all the conditions simultaneously.

By the Cohen structure theorem, fix a coefficient field
$\Bbbk\subseteq B$. For a nonzero $C_i$ with $i>0$, the parameter
ideal $Q_i$ makes $C_i$ finite over
$T_i=\Bbbk[[Y_1,\ldots,Y_i]]$, acting by $Y_j\mapsto y_j$.
Indeed, since $Q_i$ is a parameter ideal on $C_i$, the quotient
$C_i/Q_iC_i$ has finite length and hence is finite-dimensional
over the coefficient field $\Bbbk$. Moreover,
$\m^{N_i}C_i\subseteq Q_iC_i$ for some $N_i\ge1$, so the
$Q_i$-adic and $\m$-adic topologies on $C_i$ agree. Choose
$u_1,\ldots,u_t\in C_i$ whose images form a $\Bbbk$-basis of
$C_i/Q_iC_i$. Every $c\in C_i$ can be written as
\[
 c=\sum_{h=1}^{t}a_hu_h+\sum_{j=1}^{i}y_jc_j,
 \qquad a_h\in\Bbbk,\quad c_j\in C_i.
\]
Repeating this decomposition on the $c_j$ produces terms of
increasing total degree in the $y_j$, with the remainder after
$n$ steps lying in $Q_i^nC_i$. Since $C_i$ is $Q_i$-adically complete,
passing to the limit gives
\[
 c=\sum_{h=1}^{t}f_h(y_1,\ldots,y_i)u_h,
 \qquad f_h\in T_i.
\]
Thus $C_i=\sum_{h=1}^{t}T_i u_h$, proving that $C_i$ is finite
over $T_i$.

Put $A_i=B/\Ann_B(C_i)$. Multiplication on $C_i$ embeds $A_i$
into the finite $T_i$-module $\operatorname{End}_{T_i}(C_i)$,
so $A_i$ is finite over $T_i$. Define $\psi_i:T_i\to A_i$ by
$\psi_i(f)=f(y_1,\ldots,y_i)+\Ann_B(C_i)$.
The finite extension $T_i/\ker\psi_i\subseteq A_i$ preserves Krull
dimension \cite[Section~8.1, Axiom~D3]{Eisenbud1995}.
Thus $\dim(T_i/\ker\psi_i)=\dim A_i=\dim_B C_i=i$.
Since $T_i$ is a domain of dimension $i$, a nonzero kernel would
make this dimension less than $i$. Thus $\psi_i$ is injective,
and evaluation $Y_j\mapsto y_j$ identifies $T_i$ with the
regular local subring
\[
 S_i=\Bbbk[[y_1,\ldots,y_i]]\subseteq B,
 \qquad \mathfrak n_i=(y_1,\ldots,y_i)S_i.
\]
Let $\mathfrak p\in\Ass_B C_i$ and put
$\mathfrak q=\mathfrak p\cap S_i$. Since $A_i$ is finite over
$S_i$, the extension $S_i/\mathfrak q\subseteq B/\mathfrak p$
is finite. Hence
$\dim(S_i/\mathfrak q)=\dim(B/\mathfrak p)=i$.
As $S_i$ is a domain of dimension $i$, this forces
$\mathfrak q=0$. Thus every nonzero element of $S_i$ avoids
all associated primes of $C_i$ and is therefore a
nonzerodivisor on $C_i$. Consequently, $C_i$ is torsion-free
over $S_i$.
Lengths over $S_i$ and $B$ agree for the quotients by powers of
$Q_i$, since the residue field is the same. The rank formula for
multiplicity and the reduction equality therefore give
\[
 \operatorname{rank}_{S_i}C_i
 =e_{\mathfrak n_i}(C_i)=e_{Q_i}(C_i)=e_{\m}(C_i)=e_i.
\]

Choose $r_i\ge0$ with
$\m^{r_i+1}C_i=Q_i\m^{r_i}C_i$, and set $L_i=\m^{r_i}C_i$.
Since $C_i$ is torsion-free over $S_i$, multiplication by
$y_1^{r_i}$ gives $C_i\cong y_1^{r_i}C_i$. The inclusions
$y_1^{r_i}C_i\subseteq L_i\subseteq C_i$ then imply
$\operatorname{rank}_{S_i}L_i=e_i$.
For any $z\in\m$, multiplication by $z$ satisfies
\[
 zL_i\subseteq\m^{r_i+1}C_i=Q_iL_i=\mathfrak n_iL_i.
\]
Lemma~\ref{lem:characteristic-coefficients} gives a monic polynomial
$F_{i,z}\in S_i[Z]$ of degree $e_i$ such that
$F_{i,z}(z)C_i=0$, with its coefficient of $Z^{e_i-j}$ in
$\mathfrak n_i^j\subseteq Q^j$.

For $i=0$, set $F_{0,z}(Z)=Z^{e_0}$, and set $F_{i,z}=1$ for
every zero factor $C_i$. Each $F_{i,z}(z)$ maps $D_i$ into
$D_{i-1}$. Applying these maps successively from $D_d$ down to
$D_{-1}=0$ shows that
\[
 P_z(Z)=\prod_{i=0}^d F_{i,z}(Z)
\]
annihilates $M$ after evaluation at $z$.
Its degree is $\sum_i e_i=\adeg(M)$, independent of $z$.
The coefficient of $Z^{a-j}$ in $P_z$ is a sum of products
of coefficients of $Z^{e_i-j_i}$ in $F_{i,z}$, with
$0\le j_i\le e_i$ and $\sum_i j_i=j$; each product
lies in $\prod_i Q^{j_i}=Q^j$.
\end{proof}

\begin{corollary}
\label{cor:prime-power-annihilator}
Suppose that $\operatorname{char}(B/\m)=p>0$, and let $M$ be a
nonzero finite $B$-module of dimension $d\ge1$ such that
$p^vM=0$ for some $v\ge1$.
For general $y_1,\ldots,y_d\in\m$, put $Q=(y_1,\ldots,y_d)$ and
$a=\adeg(M)$. For every $z\in\m$, there is a monic polynomial
\[
 P_z(Z)=Z^a+\alpha_{1,z}Z^{a-1}+\cdots+\alpha_{a,z}\in B[Z],
 \qquad \alpha_{j,z}\in Q^j,
 \qquad P_z(z)M=0.
\]
For finite-length and zero modules, use the choices in
Proposition~\ref{prop:equicharacteristic-annihilator}.
\end{corollary}

\begin{proof}
Let $C_i=D_i/D_{i-1}$ be the factors of the dimension filtration of
$M$. For each nonzero $C_i$, set
\[
 T_{i,r}=0:_{C_i}p^r\quad(0\le r\le v),\qquad
 L_{i,r}=T_{i,r}/T_{i,r-1}\quad(1\le r\le v).
\]
Then $T_{i,0}=0$, $T_{i,v}=C_i$, and $pL_{i,r}=0$.
Multiplication by $p^{r-1}$ induces an injection
\[
 L_{i,r}\longrightarrow C_i,\qquad
 u+T_{i,r-1}\longmapsto p^{r-1}u.
\]
Thus $\Ass_B L_{i,r}\subseteq\Ass_B C_i$, so every nonzero
$L_{i,r}$ has all its associated primes of dimension $i$.
Consequently $\adeg(L_{i,r})=e(L_{i,r})$, and additivity of
multiplicity along this filtration gives
\[
 \sum_{\substack{1\le r\le v\\L_{i,r}\ne0}}e(L_{i,r})=e(C_i).
\]
In dimension zero, all these multiplicities mean lengths.

Each $L_{i,r}$ is a module over the complete equicharacteristic
ring $B/pB$. Its lengths, dimension, and multiplicity are the same
over $B$ and $B/pB$. For $i>0$, apply
Proposition~\ref{prop:equicharacteristic-annihilator} over $B/pB$
to every nonzero $L_{i,r}$, using the images of
$y_1,\ldots,y_i$. All required general choices can be imposed
simultaneously: the map
$\m/\m^2\to(\m/pB)/(\m/pB)^2$ is surjective, and there are
only finitely many layers.
Lift each coefficient from the corresponding power of the image
of $Q_i=(y_1,\ldots,y_i)B$ to that power of $Q_i$.
For each $z\in\m$, this gives a monic polynomial $F_{i,r,z}\in B[Z]$
of degree $e(L_{i,r})$ whose evaluation annihilates $L_{i,r}$,
with its coefficient of degree $e(L_{i,r})-l$ in $Q_i^l$
for $1\le l\le e(L_{i,r})$.
For $i=0$, take $F_{0,r,z}=Z^{\ell(L_{0,r})}$, and use $1$
for zero layers.

The product $F_{i,z}=\prod_{r=1}^vF_{i,r,z}$ annihilates $C_i$
after evaluation: apply $F_{i,v,z}(z),\ldots,F_{i,1,z}(z)$
successively along the filtration $T_{i,r}$.
Its degree is $e(C_i)$.
Multiplying these polynomials over all nonzero $C_i$ gives
a polynomial annihilating $M$ of degree
$\sum_{C_i\ne0}e(C_i)=\adeg(M)$. Since $Q_i\subseteq Q$, the coefficient
of degree $\adeg(M)-l$ in this product belongs to $Q^l$.
\end{proof}

\section{Triangular parameter ideals}\label{sec:parameters}

We next turn annihilating equations for the deficiency modules into a
parameter ideal whose associated graded ring has regularity zero.  The
construction is triangular: the equation used for the $i$th parameter has
coefficients in the ideal generated by the preceding chosen $y$-parameters. This
feature will also give an identity between powers in the ambient ring.

We use the associated graded notation introduced in
Section~\ref{sec:module-annihilators}. For a finite graded module $L$
over a standard graded ring $G$, its regularity is
\[
 \operatorname{reg}L
 =\max_{i\ge0}\{\operatorname{end}H^i_{G_+}(L)+i\},
\]
where local cohomology is taken with support in the positive-degree
ideal $G_+$. For a graded module $V$, we set
$\operatorname{end}V=\sup\{n:V_n\ne0\}$, with
$\operatorname{end}(0)=-\infty$.
Following \cite[Section~5]{Trung1998}, a sequence
$z_1,\ldots,z_s$ in a local ring $B$ is a \emph{$d$-sequence} if, for
$\mathfrak q=(z_1,\ldots,z_s)B$,
\begin{equation}\label{eq:dsequence-definition}
 \bigl((z_1,\ldots,z_{i-1}):z_i\bigr)\cap\mathfrak q
 =(z_1,\ldots,z_{i-1})\qquad(1\leq i\leq s).
\end{equation}
All ideals and colons in this definition are taken in $B$.
Minimality of the generating sequence is not part of this definition.

\begin{lemma}\label{lem:ann-propagation}
Let $(A,\m)$ be a local ring of dimension $s$, and put
$\mathfrak a_j=\operatorname{Ann}_A H^j_{\m}(A)$ for $0\le j<s$.
If $z_1,\ldots,z_l\in\m$ is filter regular on $A$, with $0\le l\le s$,
then
\[
 \prod_{j=0}^{l}\mathfrak a_{h+j}^{\binom lj}
 \subseteq
 \operatorname{Ann}_A H^h_{\m}\bigl(A/(z_1,\ldots,z_l)\bigr)
 \qquad(0\le h<s-l).
\]
\end{lemma}
\begin{proof}
Induct on $l$, with $l=0$ immediate. Put
$M=A/(z_1,\ldots,z_{l-1})$ and $z=z_l$.
The module $L=0:_Mz$ has finite length because $z$ is filter regular.
Apply local cohomology to the exact sequence
\[
 0\longrightarrow M/L\xrightarrow{z}M
 \longrightarrow M/zM\longrightarrow0.
\]
Since $H^0_{\m}(M/L)=H^0_{\m}(M)/L$ and
$H^j_{\m}(M/L)\cong H^j_{\m}(M)$ for $j>0$, we obtain
\begin{equation}\label{eq:lc-filter-section}
0\longrightarrow H^h_{\m}(M)/zH^h_{\m}(M)
 \longrightarrow H^h_{\m}(M/zM)
 \longrightarrow 0:_{H^{h+1}_{\m}(M)}z\longrightarrow0.
\end{equation}
This sequence also holds for $h=0$.
For an exact sequence $0\to U\to V\to W\to0$, one has
$(\operatorname{Ann}U)(\operatorname{Ann}W)\subseteq
\operatorname{Ann}V$.  Apply this observation to
\eqref{eq:lc-filter-section}. The induction hypothesis for $l-1$,
applied in cohomological degrees $h$ and $h+1$, gives exponents
$\binom{l-1}{j}$ and
$\binom{l-1}{j-1}$ for $\mathfrak a_{h+j}$. Their sum is
$\binom lj$, with the convention
$\binom{l-1}{-1}=\binom{l-1}{l}=0$.
This proves the assertion.
\end{proof}

\begin{lemma}\label{lem:torsion-dsequence}
Let $(A,\m)$ be a local ring, and let $z_1,\ldots,z_r\in\m$
be a filter regular sequence. Assume that
\begin{equation}\label{eq:successive-torsion-killing}
 z_i H^0_{\m}\bigl(A/(z_1,\ldots,z_{i-1})\bigr)=0
 \qquad(1\le i\le r).
\end{equation}
Then $z_1,\ldots,z_r$ is a $d$-sequence, and
$\mathfrak q=(z_1,\ldots,z_r)A$ satisfies
$\operatorname{reg}\operatorname{gr}_{\mathfrak q}(A)=0$.
\end{lemma}
\begin{proof}
For $1\le i\le r+1$, write $B_i=A/(z_1,\ldots,z_{i-1})$ and
$T_i=H^0_{\m}(B_i)$.
By \eqref{eq:filter-regular-torsion}, $z_i$ is regular on $B_i/T_i$
and $T_i\cap z_iB_i=z_iT_i$. The hypothesis $z_iT_i=0$ gives
\[
 0:_{B_i}z_i=T_i,\qquad T_i\cap z_iB_i=0.
\]
Consequently, the natural quotient map
$B_i\to B_{i+1}=B_i/z_iB_i$ restricts to an injection
$T_i\hookrightarrow T_{i+1}$.
Composing these quotient maps shows that
$B_i\to B_{r+1}=A/\mathfrak q$ is injective on $T_i$.
Since the kernel of this map is $\mathfrak qB_i$, the
injectivity of its restriction to $T_i$ gives
\[
 (0:_{B_i}z_i)\cap\mathfrak qB_i
 =T_i\cap\mathfrak qB_i=0.
\]
Lifting this equality to $A$ gives
\eqref{eq:dsequence-definition}.

The regularity assertion is
\cite[Corollary~5.2(i)]{Trung1998}, which applies to any ideal
generated by a $d$-sequence and imposes no restriction on the
characteristic or on the residue field.
\end{proof}

\begin{proposition}\label{prop:triangular-parameters}
Let $(R,\m)$ be a local ring, let $I$ be a proper ideal of $R$, and
put $A=R/I$, with $s=\dim A>0$. Let $y_1,\ldots,y_s\in\m$ have images forming a filter regular
system of parameters of $A$.  Write $Q_j=(y_1,\ldots,y_j)R$, with
$Q_0=0$, and $Q=Q_s$.
Suppose that there are integers $a_0,\ldots,a_{s-1}\ge0$ and, for
$0\le j<i\le s$, monic polynomials $P_{j,i}(Z)\in R[Z]$
\begin{equation}\label{eq:triangular-annihilators}
 P_{j,i}(Z)=Z^{a_j}+c_{j,i,1}Z^{a_j-1}+\cdots+c_{j,i,a_j},
 \qquad c_{j,i,l}\in Q_j^l\quad(1\le l\le a_j),
\end{equation}
such that $P_{j,i}(y_i)$ annihilates $H^j_{\m}(A)$.  If $a_j=0$,
interpret $P_{j,i}=1$, so that the hypothesis requires
$H^j_{\m}(A)=0$.
For $1\le i\le s$, set
\[
 b_i=\sum_{j=0}^{i-1}\binom{i-1}{j}a_j,
 \qquad t=\max\{1,b_s\},
 \qquad
 z_i=y_i^{\,t-b_i}
       \prod_{j=0}^{i-1}P_{j,i}(y_i)^{\binom{i-1}{j}},
 \qquad \mathfrak q=(z_1,\ldots,z_s)R.
\]
Then the following hold:
\begin{enumerate}
\item The images of $z_1,\ldots,z_s$ form a filter regular system of
parameters and a $d$-sequence in $A$.
\item $\operatorname{reg}\operatorname{gr}_{\mathfrak q}(A)=0$.
\item For each $i$, one has $z_i=F_i(y_1,\ldots,y_i)$ for a polynomial
$F_i\in R[Y_1,\ldots,Y_i]$ that is homogeneous of degree $t$
and monic in $Y_i$, where the $Y_j$ are independent indeterminates.
In the ambient ring $R$,
\begin{equation}\label{eq:ambient-balanced-identity}
 Q^{s(t-1)+1}=\mathfrak q Q^{(s-1)(t-1)}.
\end{equation}
\end{enumerate}
\end{proposition}
\begin{proof}
Since $\binom{i-1}{j}\le\binom{s-1}{j}$ for $j<i$, we have $b_i\le b_s$;
all exponents in the definition of $z_i$ are therefore nonnegative.
The coefficient condition in \eqref{eq:triangular-annihilators} implies
\begin{equation}\label{eq:triangular-prefix}
 z_i\equiv y_i^t\pmod{Q_{i-1}}.
\end{equation}
Induction on $i$ gives
\begin{equation}\label{eq:prefix-radicals}
 \sqrt{(z_1,\ldots,z_i)A}=\sqrt{(y_1,\ldots,y_i)A}
 \qquad(1\le i\le s).
\end{equation}
In particular, the full sequence is a system of parameters.

We now verify that $z_1,\ldots,z_s$ is filter regular on $A$.
Let $\mathfrak p\ne\m A$ be a prime of $A$
containing $(z_1,\ldots,z_i)A$.  By \eqref{eq:prefix-radicals}, it
contains $(y_1,\ldots,y_i)A$.  Filter regularity of the latter sequence
shows that $y_1,\ldots,y_i$ is a regular sequence in $A_{\mathfrak p}$.
The ideal $(y_1,\ldots,y_i)A_{\mathfrak p}$ thus has grade $i$.
Localizing \eqref{eq:prefix-radicals} shows that
$(z_1,\ldots,z_i)A_{\mathfrak p}$ has the same radical as
$(y_1,\ldots,y_i)A_{\mathfrak p}$. Since ideals with the same radical have the same grade,
$(z_1,\ldots,z_i)A_{\mathfrak p}$ also has grade $i$.
In a Noetherian local ring, a sequence of $i$ elements generating an
ideal of grade $i$ is regular
\cite[Theorems~16.5(ii) and~16.8]{Matsumura1989}. Hence the ordered
sequence $z_1,\ldots,z_i$ is regular in $A_{\mathfrak p}$.
It follows that $z_i$ avoids every nonmaximal associated prime of
$A/(z_1,\ldots,z_{i-1})A$: a prime containing $z_i$ would contradict
this regularity.
This proves filter regularity of the full sequence.

Put $\mathfrak a_j=\operatorname{Ann}_A H^j_{\m}(A)$.
The annihilation assumptions and the product defining $z_i$ give
\[
 z_i\in\prod_{j=0}^{i-1}\mathfrak a_j^{\binom{i-1}{j}}
 \quad\text{in }A.
\]
Lemma~\ref{lem:ann-propagation}, applied to the filter regular prefix
$z_1,\ldots,z_{i-1}$, now shows that
\[
 z_iH^0_{\m}\bigl(A/(z_1,\ldots,z_{i-1})A\bigr)=0.
\]
Lemma~\ref{lem:torsion-dsequence} proves the first two assertions.

We now prove (iii). Express each coefficient $c_{j,i,l}\in Q_j^l$
as a homogeneous polynomial of degree $l$ in $y_1,\ldots,y_j$, with
coefficients in $R$. Put $\mathbf Y=(Y_1,\ldots,Y_s)$, with the $Y_j$
independent indeterminates over $R$. Substituting these expressions
into the formula defining $z_i$ and replacing the $y_k$ by $Y_k$,
while keeping the coefficients in $R$ fixed, gives a polynomial
$F_i\in R[Y_1,\ldots,Y_i]$ that is homogeneous of degree $t$
and monic of degree $t$ in $Y_i$, with $F_i(y_1,\ldots,y_i)=z_i$.
Set $B_0=R$ and
\[
 B_i=R[Y_1,\ldots,Y_i]/(F_1,\ldots,F_i)
 \qquad(1\le i\le s).
\]
Since $F_1,\ldots,F_{i-1}$ involve only $Y_1,\ldots,Y_{i-1}$, we have
\[
 B_i\cong B_{i-1}[Y_i]/(\widetilde F_i),
\]
where $\widetilde F_i$ is obtained by mapping the coefficients of $F_i$,
viewed as a polynomial in $Y_i$, to $B_{i-1}$.
It remains monic of degree $t$. Division by a monic polynomial over
any coefficient ring gives a unique remainder of degree less than $t$.
Thus $B_i$ is a free $B_{i-1}$-module with basis consisting of the
residue classes of $1,Y_i,\ldots,Y_i^{t-1}$.
Induction on $i$ shows that $B_s$ is a free $R$-module with basis
consisting of the residue classes of the monomials
$Y_1^{u_1}\cdots Y_s^{u_s}$, where $0\le u_i<t$.
Since the $F_i$ are homogeneous, these basis elements are homogeneous
of degree at most $s(t-1)$. For $b=s(t-1)+1$, we thus have
$(B_s)_b=0$, and hence
\[
 R[\mathbf Y]_b=\sum_{i=1}^s F_iR[\mathbf Y]_{b-t}.
\]
Evaluation at $Y_i=y_i$ gives $Q^b=\mathfrak qQ^{b-t}$.
Since $b-t=(s-1)(t-1)$, this is precisely
\eqref{eq:ambient-balanced-identity}.
\end{proof}

\section{From parameter ideals to the maximal-ideal Hilbert function}
\label{sec:transport}

Throughout this section the local rings have infinite residue fields.
This suffices for both main theorems, since the faithfully flat
reduction to this case is made in Section~\ref{sec:main-proof}.
We prove a perturbation criterion that is independent of homological
degree. Its input consists of a reduction of the maximal ideal and a parameter
ideal whose associated graded ring has regularity zero. The construction in
the preceding section will supply these two ideals simultaneously.
The comparison with a parameter filtration follows the proofs of
\cite[Proposition~1]{ST97} and \cite[Theorem~4.2]{QVT}.
Here we construct the required parameter ideal from the final quotient
and use successive generic approximations to carry it through the proof.
The argument has three parts: a one-element perturbation lemma,
a comparison of the two filtrations using an intersection identity,
and the passage to several generators.

All regularity comparisons below may be made over a common
associated graded ring. If $\mathfrak q$ is an ideal of a local ring $W$
and $C$ is a quotient of $W$, then
$\operatorname{gr}_{\mathfrak q}(C)$ is a quotient of
$\operatorname{gr}_{\mathfrak q}(W)$. Images of the same degree-one
generators compute the local cohomology of either module by the
\v{C}ech complex. Thus regularity is unchanged by this change of
presentation, and the usual exact-sequence inequalities apply.

\begin{lemma}\label{lem:tr-finite-part}
Let $(B,\mathfrak n)$ be local, let $\mathfrak a$ be an
$\mathfrak n$-primary ideal, and assume that
$\operatorname{reg}\operatorname{gr}_{\mathfrak a}(B)=0$.
For every proper ideal $L$ of $B$ contained in $H^0_{\mathfrak n}(B)$,
one has
\[
 L\cap\mathfrak a=0,\qquad
 \operatorname{reg}\operatorname{gr}_{\mathfrak a}(B/L)=0.
\]
Moreover, there is an exact sequence of graded modules
\begin{equation}\label{eq:tr-finite-graded}
 0\longrightarrow L\longrightarrow
 \operatorname{gr}_{\mathfrak a}(B)\longrightarrow
 \operatorname{gr}_{\mathfrak a}(B/L)\longrightarrow0,
\end{equation}
where $L$ is placed in degree zero.
\end{lemma}

\begin{proof}
The induced filtration on $L$ gives a graded submodule
\[
 E=\bigoplus_{n\geq0}
 \frac{L\cap\mathfrak a^n}{L\cap\mathfrak a^{n+1}}
 \ \subseteq\ G:=\operatorname{gr}_{\mathfrak a}(B).
\]
Since a power of $\mathfrak a$ annihilates $L$, a power of $G_+$ annihilates
$E$. Thus $E\subseteq H^0_{G_+}(G)$. Regularity zero implies $E_n=0$ for
$n\geq1$, and hence
$L\cap\mathfrak a=L\cap\mathfrak a^n$ for all $n\geq1$.
Krull intersection gives $L\cap\mathfrak a=0$.
The kernel of the natural graded surjection is therefore
$L$ in degree zero, proving \eqref{eq:tr-finite-graded}.
The module $L$ is concentrated in degree zero and hence killed by $G_+$,
so it has regularity zero if nonzero (and regularity $-\infty$
if zero). The exact-sequence inequality therefore gives
\[
 \operatorname{reg}\operatorname{gr}_{\mathfrak a}(B/L)
 \leq\max\{\operatorname{reg}G,\operatorname{reg}L-1\}\leq0.
\]
This is a nonzero standard graded ring, so its regularity is
nonnegative by \cite[Corollary~2.3(ii)]{Trung1998}.
The asserted equality follows.
\end{proof}

\begin{lemma}\label{lem:tr-one-element}
Let $(B,\mathfrak n)$ be a local ring and let $x\in\mathfrak n$ be filter
regular. Put $C=B/xB$ and assume $\dim C>0$. Suppose that
$\mathfrak q\subseteq\mathfrak n$ is an ideal whose image in $C$
is $\mathfrak nC$-primary and that
\[
 \operatorname{reg}\operatorname{gr}_{\mathfrak q}(C)=0.
\]
For every $\eta\in\mathfrak q^2$, the element $x'=x+\eta$ is filter
regular, and $C'=B/x'B$ satisfies
\begin{align}
 \operatorname{reg}\operatorname{gr}_{\mathfrak q}(C')&=0,
 \label{eq:tr-principal-reg}\\
 \ell_B(C'/\mathfrak q^nC')&=\ell_B(C/\mathfrak q^nC)
 \quad(n\geq1).
 \label{eq:tr-principal-hf}
\end{align}
\end{lemma}

\begin{proof}
Set $T=H^0_{\mathfrak n}(B)$ and $\overline B=B/T$.
The torsion observation \eqref{eq:filter-regular-torsion} shows
that $x$ is regular on $\overline B$ and $T\cap xB=xT$, so
\[
 L=T/xT\ \subseteq\ C,
 \qquad \overline C=C/L=\overline B/x\overline B.
\]
The ideal $L$ is proper because it has finite length and $\dim C>0$.
Lemma~\ref{lem:tr-finite-part} gives
\begin{equation}\label{eq:tr-torsion-data}
 \mathfrak qL=0,\qquad L\cap\mathfrak q C=0,
 \qquad\operatorname{reg}\operatorname{gr}_{\mathfrak q}(\overline C)=0.
\end{equation}
In particular $\mathfrak qT\subseteq xT$. It follows that
\[
 \eta T\subseteq\mathfrak q^2T
 \subseteq x\mathfrak qT\subseteq\mathfrak n xT.
\]
Thus $x'T\subseteq xT$ and $xT\subseteq x'T+\mathfrak n xT$.
Nakayama's lemma gives
\begin{equation}\label{eq:tr-fixed-torsion-image}
 x'T=xT.
\end{equation}

We next compare the quotients of $\overline B$. For an ideal $J$ in a
local ring $V$ and an ideal $\mathfrak a$ of $V$, write
$\operatorname{ar}_{\mathfrak a}(J)$ for the Artin--Rees number of $J$
with respect to $\mathfrak a$: the least integer $c\geq0$ such that
$J\cap\mathfrak a^{n+1}=\mathfrak a(J\cap\mathfrak a^n)$ for $n\geq c$.
The regularity bound
\cite[Corollary 2.5]{QT} gives
\[
 \operatorname{ar}_{\mathfrak q}(x\overline B)
 \leq\operatorname{reg}\operatorname{gr}_{\mathfrak q}
            (\overline B/x\overline B)+1\leq1.
\]
Since $x$ is regular on $\overline B$ and
$\eta\in\mathfrak q^2$, the Artin--Rees bound above allows us
to apply \cite[Proposition~3.2]{QT} over $\overline B$
with respect to $\mathfrak q\overline B$.
It follows that $x'$ is regular on $\overline B$ and that
$(x)$ and $(x')$ have the same initial ideal in
$\operatorname{gr}_{\mathfrak q}(\overline B)$.
Taking the corresponding graded quotients gives
\begin{equation}\label{eq:tr-reduced-graded}
 \operatorname{gr}_{\mathfrak q}(\overline B/x'\overline B)
 \cong\operatorname{gr}_{\mathfrak q}(\overline B/x\overline B).
\end{equation}

Since $x'$ is regular modulo $T$, its kernel on $B$ is contained in $T$;
thus $x'$ is filter regular. Moreover,
$T\cap x'B=x'T=xT$. We may therefore identify
$L'=T/x'T$ with the same finite module $L=T/xT$.
Because $\eta\in\mathfrak q$, the ideals
$x'B+\mathfrak q$ and $xB+\mathfrak q$ coincide. In particular,
$\mathfrak qC'$ is primary for the maximal ideal of $C'$.
The middle assertion
of \eqref{eq:tr-torsion-data} says exactly that
\[
 T\cap(xB+\mathfrak q)=xT.
\]
Using \eqref{eq:tr-fixed-torsion-image}, we obtain
$L'\cap\mathfrak qC'=0$ and $\mathfrak qL'=0$.
Hence there is an exact sequence
\[
 0\longrightarrow L\longrightarrow
 \operatorname{gr}_{\mathfrak q}(C')\longrightarrow
 \operatorname{gr}_{\mathfrak q}(\overline B/x'\overline B)
 \longrightarrow0,
\]
again with $L$ in degree zero. The regularity inequality for this
exact sequence and \eqref{eq:tr-reduced-graded} give
$\operatorname{reg}\operatorname{gr}_{\mathfrak q}(C')\leq0$.
It is exactly zero because $C'\neq0$ and a nonzero
standard graded ring has nonnegative regularity, as used in
Lemma~\ref{lem:tr-finite-part}. This proves
\eqref{eq:tr-principal-reg}.

For every $n\geq1$, reduce the exact sequences
\[
 0\to L\to C\to\overline C\to0,
 \qquad
 0\to L'\to C'\to\overline B/x'\overline B\to0
\]
modulo $\mathfrak q^n$.
Since $L\cap\mathfrak q^nC=0$ and
$L'\cap\mathfrak q^nC'=0$, both resulting sequences remain
exact. Taking lengths and using $L'\cong L$ gives
\[
 \ell(C/\mathfrak q^nC)
 =\ell(L)+\ell(\overline C/\mathfrak q^n\overline C),
\]
and the analogous equality for $C'$. By \eqref{eq:tr-reduced-graded},
the quotients $\overline B/x'\overline B$ and $\overline C$
have the same $\mathfrak q$-adic Hilbert function.
The two length formulas above, with the common term $\ell(L)$,
now give \eqref{eq:tr-principal-hf}.
\end{proof}

\begin{lemma}
\label{lem:tr-parameter-intersection}
Let $(R,\mathfrak m)$ be local with infinite residue field, let $J$
be a proper ideal, and put $B=R/J$. Suppose that $s=\dim B>0$, that
$\mathfrak q=(z_1,\ldots,z_s)\subseteq\mathfrak m$ is an ideal
whose image in $B$ is $\mathfrak mB$-primary, and that
$\operatorname{reg}\operatorname{gr}_{\mathfrak q}(B)=0$.
Then, for $\mathfrak b=J+\mathfrak q$,
\begin{equation}\label{eq:tr-intersection}
 \mathfrak b^{v+1}\cap J=J\mathfrak b^v\qquad(v\geq0).
\end{equation}
\end{lemma}

\begin{proof}
The ideal $\mathfrak qB$ requires exactly $s$ generators: it is
generated by the images of $z_1,\ldots,z_s$ in $B$, and an ideal primary
for the maximal ideal of a ring of dimension $s$ cannot be generated
by fewer, by the height theorem. By \cite[Corollary~5.2(ii)]{Trung1998}, it has a generating
$d$-sequence $\bar w_1,\ldots,\bar w_s$. Here the generators may
be taken minimal: the proof of that corollary uses the minimal
basis supplied by \cite[Theorem~4.9]{Trung1998}.
Let $\bar z_i$ be the image of $z_i$ in $B$. Expressing the new
generators in terms of the old ones gives
$(\bar w_i)_i=U(\bar z_i)_i$, where $U$ is an $s\times s$ matrix over
$B$ and both tuples are viewed as column vectors. The images of both
generating sets are bases of
$\mathfrak qB/(\mathfrak mB)\mathfrak qB$. The reduction of
$U$ modulo $\mathfrak mB$ is therefore invertible, so $U$
itself is invertible over the local ring $B$.
Lift its entries to a matrix $\widetilde U$ over $R$. Its
determinant remains a unit, since its residue modulo
$\mathfrak m$ is nonzero. Define lifts $w_1,\ldots,w_s\in R$ of
$\bar w_1,\ldots,\bar w_s$ by setting
$(w_i)_i=\widetilde U(z_i)_i$.
Then $(w_1,\ldots,w_s)R=\mathfrak q$, and the images of these
elements form a $d$-sequence modulo $J$.

We prove $\mathfrak q^n\cap J\subseteq J\mathfrak q^{n-1}$
for $n\geq1$ directly from \eqref{eq:dsequence-definition},
following the proof of \cite[Theorem~2.1]{Huneke82}.
More generally, this inclusion holds whenever the images of a
sequence generating $\mathfrak q$ are a $d$-sequence modulo $J$.
Induct on the number of generators, the empty case being trivial;
for a fixed sequence, induct also on $n$, with $n=1$ immediate.
Write $x=w_1$ and $\mathfrak a=(w_2,\ldots,w_s)$.
The tail is a $d$-sequence modulo $J+(x)$ by the same definition,
so the induction on the number of generators gives
\[
 \mathfrak a^n\cap(J+(x))\subseteq(J+(x))\mathfrak a^{n-1}.
\]
For $a\in\mathfrak q^n\cap J$, with $n\geq2$, write
$a=b+xc$, where $b\in\mathfrak a^n$ and $c\in\mathfrak q^{n-1}$.
Then $b\in J+(x)$, and the displayed inclusion gives
$b=u+xv$, with $u\in J\mathfrak a^{n-1}$ and
$v\in\mathfrak a^{n-1}$. Hence
$x(c+v)=a-u\in J$. The first condition in
\eqref{eq:dsequence-definition}, lifted modulo $J$, is
\[
 (J:x)\cap(J+\mathfrak q)=J.
\]
Since $c+v\in\mathfrak q^{n-1}\subseteq\mathfrak q$, it follows
that $c+v\in J\cap\mathfrak q^{n-1}$.
Induction on $n$ now gives $c+v\in J\mathfrak q^{n-2}$, and thus
$a=u+x(c+v)\in J\mathfrak q^{n-1}$, as required.

Finally, expand $\mathfrak b^{v+1}=J\mathfrak b^v+\mathfrak q^{v+1}$. If an element
of this sum lies in $J$, its $\mathfrak q^{v+1}$ part lies in
$J\cap\mathfrak q^{v+1}\subseteq J\mathfrak q^v\subseteq J\mathfrak b^v$.
This proves $\mathfrak b^{v+1}\cap J\subseteq J\mathfrak b^v$; the reverse inclusion
follows from $J\subseteq \mathfrak b$.
\end{proof}

When applying this lemma to two quotients $R/J$ and $R/J'$, we may
choose generating sequences
\[
 \mathfrak q=(w_1,\ldots,w_s)R=(w'_1,\ldots,w'_s)R
\]
whose images are $d$-sequences in $R/J$ and $R/J'$, respectively.
These choices may be made independently, since the ideal
$\mathfrak q$, its powers, and all identities involving these powers
remain unchanged.

\begin{proposition}\label{prop:tr-hilbert-comparison}
Let $(R,\mathfrak m)$ be local and let $J,J'$ be proper ideals.
Put $s=\dim R/J>0$, let
$Q=(y_1,\ldots,y_s)\subseteq\mathfrak m$ and let $t\geq1$ be an integer,
and put $b=s(t-1)+1$. Suppose that an ideal
$\mathfrak q=(z_1,\ldots,z_s)\subseteq Q^t$
satisfies the identity
\begin{equation}\label{eq:tr-balanced-identity}
 Q^b=\mathfrak qQ^{b-t}
\end{equation}
in $R$. Assume that, for an integer $\rho\geq0$,
\begin{equation}\label{eq:tr-original-reduction}
 \mathfrak m^{\rho+1}(R/J)=Q\mathfrak m^\rho(R/J).
\end{equation}
Set $N=\rho+b+t+1$, and suppose that
\begin{enumerate}
 \item $J+\mathfrak m^N=J'+\mathfrak m^N$;
 \item $J+\mathfrak q$ and $J'+\mathfrak q$ are
 $\mathfrak m$-primary, and
 \[
 \ell\bigl(R/(J+\mathfrak q^n)\bigr)
 =\ell\bigl(R/(J'+\mathfrak q^n)\bigr)\quad(n\geq1);
 \]
 \item $\operatorname{reg}\operatorname{gr}_{\mathfrak q}(R/J)
       =\operatorname{reg}\operatorname{gr}_{\mathfrak q}(R/J')=0$.
\end{enumerate}
Then $R/J$ and $R/J'$ have the same $\mathfrak m$-adic Hilbert function.
\end{proposition}

\begin{proof}
The Hilbert--Samuel polynomial of a primary ideal has degree
equal to the dimension of the ring. Assumption (ii) therefore
gives $\dim R/J'=\dim R/J=s$. Since $\mathfrak q$ is generated
by $s$ elements, assumption (ii) shows that its images in $R/J$
and $R/J'$ are parameter ideals. Together with assumption (iii),
this allows us to apply Lemma~\ref{lem:tr-parameter-intersection}
to both quotients.

By \eqref{eq:tr-original-reduction} and assumption (i), we have
\[
 \mathfrak m^{\rho+2}
 \subseteq Q\mathfrak m^{\rho+1}+\mathfrak mJ
 \subseteq Q\mathfrak m^{\rho+1}+J'+\mathfrak m^{N+1}.
\]
Since $N\geq\rho+2$, we have
$\mathfrak m^{N+1}\subseteq\mathfrak m^{\rho+3}$.
Applying Nakayama's lemma to the finitely generated $R$-module
$\mathfrak m^{\rho+2}(R/J')/Q\mathfrak m^{\rho+1}(R/J')$,
we obtain
\begin{equation}\label{eq:tr-common-reduction}
 \mathfrak m^{\rho+2}(R/J')=Q\mathfrak m^{\rho+1}(R/J').
\end{equation}
On the other hand, \eqref{eq:tr-original-reduction} gives
$\mathfrak m^{\rho+b}(R/J)=Q^b\mathfrak m^\rho(R/J)$.
Using \eqref{eq:tr-balanced-identity}, we obtain
\[
 \mathfrak m^{\rho+b}\subseteq J+\mathfrak q.
\]
Hence $\mathfrak m^N\subseteq\mathfrak m(J+\mathfrak q)$.
By assumption (i) and $\mathfrak m^N\subseteq\mathfrak m(J+\mathfrak q)$,
we have $J'+\mathfrak q\subseteq J+\mathfrak q$ and
\[
 J+\mathfrak q=(J'+\mathfrak q)+\mathfrak m(J+\mathfrak q).
\]
Applying Nakayama's lemma to $(J+\mathfrak q)/(J'+\mathfrak q)$,
we obtain
\begin{equation}\label{eq:tr-common-K}
 \mathfrak b:=J+\mathfrak q=J'+\mathfrak q.
\end{equation}

For $0\leq j<t$, define
\[
 n_j=\rho+b+1+j.
\]
Then $n_j\leq N-1$ and
$n_j-t\geq\rho+1$, because $b-t=(s-1)(t-1)\geq0$.
Repeated use of \eqref{eq:tr-balanced-identity} gives
\[
 Q^{b+j+vt}=\mathfrak q^{v+1}Q^{b-t+j}
 \quad(j,v\geq0).
\]
Together with \eqref{eq:tr-original-reduction} and
\eqref{eq:tr-common-reduction}, this gives
\begin{equation}\label{eq:tr-all-residues}
 \mathfrak m^{n_j+vt}+J
   =\mathfrak q^{v+1}\mathfrak m^{n_j-t}+J,
 \qquad
 \mathfrak m^{n_j+vt}+J'
   =\mathfrak q^{v+1}\mathfrak m^{n_j-t}+J'.
\end{equation}
Here and below $0\leq j<t$ and $v\geq0$.
Indeed, the calculation in $R/J$ is
\[
 \begin{aligned}
 \mathfrak m^{n_j+vt}(R/J)
 &=Q^{b+j+vt}\mathfrak m^{\rho+1}(R/J)\\
 &=\mathfrak q^{v+1}Q^{b-t+j}\mathfrak m^{\rho+1}(R/J)\\
 &=\mathfrak q^{v+1}\mathfrak m^{n_j-t}(R/J).
 \end{aligned}
\]
The same calculation applies with $J$ replaced by $J'$.
All exponents of $Q$ in this calculation are nonnegative.

Lemma~\ref{lem:tr-parameter-intersection} yields
$\mathfrak b^{v+1}\cap J=J\mathfrak b^v$ and $\mathfrak b^{v+1}\cap J'=J'\mathfrak b^v$.
Using \eqref{eq:tr-all-residues} with $v=0$, we obtain
\begin{align*}
 \mathfrak b^{v+1}\mathfrak m^{n_j-t}+(\mathfrak b^{v+1}\cap J)
 &=\mathfrak b^v(\mathfrak b\mathfrak m^{n_j-t}+J)\\
 &=\mathfrak b^v(\mathfrak m^{n_j}+J)\\
 &=\mathfrak b^v(\mathfrak m^{n_j}+J')\\
 &=\mathfrak b^{v+1}\mathfrak m^{n_j-t}+(\mathfrak b^{v+1}\cap J').
\end{align*}
The middle equality follows from assumption (i), since $n_j\leq N$.
Put $l_j=n_j-t$ and denote the common ideal above by
\[
 \mathfrak d_{j,v}=\mathfrak b^{v+1}\mathfrak m^{l_j}+J\mathfrak b^v
        =\mathfrak b^{v+1}\mathfrak m^{l_j}+J'\mathfrak b^v
        \ \subseteq \mathfrak b^{v+1}.
\]
To compare the $\mathfrak m$-adic Hilbert functions, consider the natural surjections
\[
 \frac{R}{\mathfrak b^{v+1}\mathfrak m^{n_j-t}+J}
 \longrightarrow\frac{R}{\mathfrak b^{v+1}+J},
 \qquad
 \frac{R}{\mathfrak b^{v+1}\mathfrak m^{n_j-t}+J'}
 \longrightarrow\frac{R}{\mathfrak b^{v+1}+J'}.
\]
Both kernels are canonically identified with
$\mathfrak b^{v+1}/\mathfrak d_{j,v}$.
Modulo either $J$ or $J'$, every power of $\mathfrak b$ becomes the
corresponding power of $\mathfrak q$. The sources above have finite
length by \eqref{eq:tr-all-residues}, and therefore so do their
kernels and targets. Their targets have equal length by the assumed
equality of the $\mathfrak q$-adic Hilbert functions and \eqref{eq:tr-common-K}.
Equation \eqref{eq:tr-all-residues} therefore gives equality of the
$\mathfrak m$-adic Hilbert functions at every exponent $n_j+vt$.
These exponents comprise all integers at least $n_0$: divide
$n-n_0$ by $t$ to obtain $n=n_0+j+vt=n_j+vt$, with
$0\leq j<t$ and $v\geq0$.
For the smaller exponents the equality follows from
$J+\mathfrak m^N=J'+\mathfrak m^N$.
\end{proof}

We now use Lemma~\ref{lem:tr-one-element} to treat perturbations
of several generators. The only sequences needed are those
obtained by replacing the first
generator, then the second, and so on. We use the following
Koszul exact sequence. For a finite sequence $\mathbf v$ in a ring $B$,
write $K_\bullet(\mathbf v;B)$ for its Koszul complex and
$H_j(\mathbf v;B)$ for its $j$th homology module. For $x\in B$,
adjoining $x$ expresses the new complex as the mapping cone of
multiplication by $x$ on $K_\bullet(\mathbf v;B)$.
Its homology sequence gives
\begin{equation}\label{eq:tr-koszul-adjunction}
 0\longrightarrow
 \frac{H_j(\mathbf v;B)}{xH_j(\mathbf v;B)}
 \longrightarrow H_j(\mathbf v,x;B)
 \longrightarrow 0:_{H_{j-1}(\mathbf v;B)}x
 \longrightarrow0\qquad(j\geq1).
\end{equation}
For the empty sequence we use $H_0(\varnothing;B)=B$ and
$H_j(\varnothing;B)=0$ for $j>0$.
Permuting the entries gives an isomorphic Koszul complex, so
this sequence may be used with any one entry singled out as $x$.

\begin{lemma}
\label{lem:tr-successive-approximation}
Let $(R,\mathfrak m)$ be local with infinite residue field, and let
$f_1,\ldots,f_r\in\mathfrak m$ have finite-length positive Koszul
homology over $R$, with $r\ge1$.
For $a_1,\ldots,a_r\in\mathfrak m$ and an integer $\nu\geq1$, there are
$g_i\in a_i+\mathfrak m^\nu$ such that each sequence
\[
 (g_1,\ldots,g_i,f_{i+1},\ldots,f_r)\qquad(0\leq i\leq r)
\]
has finite-length positive Koszul homology over $R$.
\end{lemma}
\begin{proof}
Consider one replacement in a sequence $(\mathbf v,x)$ with
finite-length positive Koszul homology, and fix a target $a\in\mathfrak m$.
Write
$L_j=H_j(\mathbf v;R)$ for $j\geq0$.
By \eqref{eq:tr-koszul-adjunction}, $L_j/xL_j$ has finite length
for $j>0$. If $L_j\ne0$, this quotient is nonzero by Nakayama,
and the dimension inequality gives
\[
 \dim L_j\leq\dim(L_j/xL_j)+1\leq1\qquad(j>0).
\]

We can choose $g\in a+\mathfrak m^\nu$ filter regular on all the
finitely many nonzero $L_j$, including $L_0$.
To prove this, let $\mathcal P$ be the union of their nonmaximal
associated primes. If $\mathcal P$ is empty, take $g=a$.
Otherwise $\mathfrak m^\nu$ is contained in none of these primes,
so prime avoidance gives $u\in\mathfrak m^\nu$ outside all of them.
Choose an infinite subset $\Lambda\subseteq R$ whose elements have
pairwise distinct images in $R/\mathfrak m$.
For each $\mathfrak p\in\mathcal P$, at most one $\lambda\in\Lambda$
satisfies $a+\lambda u\in\mathfrak p$: two distinct such elements
$\lambda,\lambda'\in\Lambda$ would give
$(\lambda-\lambda')u\in\mathfrak p$, although $\lambda-\lambda'$
is a unit and $u\notin\mathfrak p$.
Since $\mathcal P$ is finite and $\Lambda$ is infinite, we may choose
$\lambda\in\Lambda$ outside these finitely many exceptions and set
$g=a+\lambda u$.
Then $0:_{L_j}g$ has finite length for every $j$.

For $j>0$, every nonmaximal prime in the support of $L_j$ is
minimal, since $\dim L_j\leq1$, and hence is associated.
Thus $g$ avoids all such primes and $L_j/gL_j$ has finite length.
In \eqref{eq:tr-koszul-adjunction} with $g$ in place of $x$,
both outer terms now have finite length, including at $j=1$.
It follows that $H_j(\mathbf v,g;R)$ has finite length for all $j>0$.

Starting with $(f_1,\ldots,f_r)$, repeat this step with
\[
 \mathbf v=(g_1,\ldots,g_{i-1},f_{i+1},\ldots,f_r),
 \qquad x=f_i,\quad a=a_i\quad(1\leq i\leq r).
\]
This constructs exactly the sequences asserted in the lemma.
\end{proof}

\begin{proposition}[Perturbation criterion]
\label{thm:tr-criterion}
Let $(R,\mathfrak m)$ be local with infinite residue field and let
$\f=(f_1,\ldots,f_r)$ be a sequence in $\mathfrak m$, with $r\ge1$,
whose positive Koszul homology modules have finite length. Put
$I=(\f)$ and $A=R/I$, and assume $s=\dim A>0$.
Suppose that $Q=(y_1,\ldots,y_s)\subseteq\mathfrak m$ satisfies
\[
 \mathfrak m^{\rho+1}A=Q\mathfrak m^\rho A
\]
for some integer $\rho\geq0$. Suppose also that, for some integer $t\geq1$, an ideal
$\mathfrak q=(z_1,\ldots,z_s)\subseteq Q^t$ satisfies
\[
 Q^{s(t-1)+1}=\mathfrak qQ^{(s-1)(t-1)}
 \quad\text{in }R,
 \qquad
 \operatorname{reg}\operatorname{gr}_{\mathfrak q}(A)=0.
\]
Put $N=\rho+(s+1)(t-1)+3$. For every
$\varepsilon_1,\ldots,\varepsilon_r\in\mathfrak m^N$, the quotients
$R/(f_1+\varepsilon_1,\ldots,f_r+\varepsilon_r)$ and $A$
have the same $\mathfrak m$-adic Hilbert function.
\end{proposition}

\begin{proof}
Put $b=s(t-1)+1$, so $N=\rho+b+t+1$. The reduction identity and the
identity for powers of $Q$ give
\begin{equation}\label{eq:tr-square-containment}
 \mathfrak m^{N-1}A
   =Q^{b+t}\mathfrak m^\rho A
   =\mathfrak q^2Q^{b-t}\mathfrak m^\rho A
   \subseteq\mathfrak q^2A.
\end{equation}
Here the second equality uses the ambient identity twice:
$Q^{b+t}=\mathfrak qQ^b=\mathfrak q^2Q^{b-t}$.
The exponent $b-t=(s-1)(t-1)$ is nonnegative, including when
$s=1$ or $t=1$.
In particular, $\mathfrak qA$ is primary for the maximal ideal
of $A$. Since it has $s=\dim A$ generators, it is a parameter
ideal. This generator count will also allow us to apply
Lemma~\ref{lem:tr-parameter-intersection} to the approximating quotients.

Fix $\varepsilon_i\in\mathfrak m^N$ and set
$f_i'=f_i+\varepsilon_i$, $I'=(f_1',\ldots,f_r')$.
For any $\nu\geq N$, apply Lemma~\ref{lem:tr-successive-approximation} with
$a_i=f_i'$ and obtain $g_i\equiv f_i'\pmod{\mathfrak m^\nu}$ in $R$.
The choice of $g_1,\ldots,g_r$ may depend on $\nu$.
For $0\leq i\leq r$ put
\[
 J_i=(g_1,\ldots,g_i,f_{i+1},\ldots,f_r)R,
 \qquad C_i=R/J_i.
\]
Since $J_i+\mathfrak m^N=I+\mathfrak m^N$,
\eqref{eq:tr-square-containment} gives
$\mathfrak m^{N-1}\subseteq J_i+\mathfrak q^2+\mathfrak m^N$.
Thus the finite module
$\mathfrak m^{N-1}(R/(J_i+\mathfrak q^2))$
equals its product with $\mathfrak m$ and is zero by Nakayama.
Consequently,
\begin{equation}\label{eq:tr-all-square-containments}
 \mathfrak m^{N-1}\subseteq J_i+\mathfrak q^2
 \qquad(0\leq i\leq r).
\end{equation}
In particular, $\mathfrak qC_i$ is primary for the maximal ideal
of $C_i$ for every $i$.

We prove by induction on $i$ that
\begin{equation}\label{eq:tr-induction-data}
 \operatorname{reg}\operatorname{gr}_{\mathfrak q}(C_i)=0,
 \qquad
 \ell(C_i/\mathfrak q^nC_i)=\ell(A/\mathfrak q^nA)
 \quad(n\geq1).
\end{equation}
These assertions hold for $i=0$ by hypothesis.
Suppose they hold for $i-1$. Equality of the Hilbert functions of
primary ideals gives $\dim C_{i-1}=\dim A=s>0$, as required in
Lemma~\ref{lem:tr-one-element}. Work in
\[
 W_i=R/(g_1,\ldots,g_{i-1},f_{i+1},\ldots,f_r)R.
\]
The sequence generating $J_{i-1}$ has finite-length positive
Koszul homology by Lemma~\ref{lem:tr-successive-approximation}.
For $x$ the image of $f_i$ in $W_i$, the case $j=1$ of
\eqref{eq:tr-koszul-adjunction} gives a surjection from this finite-length
first homology onto $0:_{W_i}x$. Thus $x$ is filter regular and
$W_i/xW_i=C_{i-1}$.
The image $\delta$ of $g_i-f_i$ belongs to $\mathfrak m^NW_i$.
Multiplying \eqref{eq:tr-all-square-containments} for $J_{i-1}$
by $\mathfrak m$ and passing to $W_i$ gives
\[
 \delta\in\mathfrak m xW_i+\mathfrak q^2W_i.
\]
Write $\delta=ax+\eta$, with $a$ in the maximal ideal and
$\eta\in\mathfrak q^2W_i$. Since $1+a$ is a unit,
\[
 (x+\delta)W_i=(x+\eta/(1+a))W_i.
\]
Lemma~\ref{lem:tr-one-element} applies to the right-hand side and
proves \eqref{eq:tr-induction-data} for $i$.
This proves the induction without increasing $N$.

We now apply Proposition~\ref{prop:tr-hilbert-comparison} over $R$
to $J_0=I$ and $J_r$. Their congruence modulo $\mathfrak m^N$
was noted above. The reduction identity for $A$ and the ambient identity for $Q$
are the hypotheses of the criterion. The parameter ideals are
primary by \eqref{eq:tr-all-square-containments}, and
\eqref{eq:tr-induction-data} supplies both equality of the
$\mathfrak q$-adic Hilbert functions and regularity zero.
Thus every hypothesis of that proposition is satisfied, and it gives
\[
 \ell_R\bigl(R/(J_r+\mathfrak m^n)\bigr)
   =\ell_R\bigl(R/(I+\mathfrak m^n)\bigr)
 \quad(n\geq1).
\]
Finally, fix $n\geq1$ and choose $\nu\geq\max\{N,n\}$.
Choose $g_1,\ldots,g_r$ as above for this value of $\nu$.
Since $g_i\equiv f_i'\pmod{\mathfrak m^\nu}$,
$J_r+\mathfrak m^n=I'+\mathfrak m^n$.
The equality already proved for $J_r$ therefore gives
\[
 \ell_R\bigl(R/(I'+\mathfrak m^n)\bigr)
 =\ell_R\bigl(R/(I+\mathfrak m^n)\bigr).
\]
The integer $n$ was arbitrary.
The choice of $g_1,\ldots,g_r$ may depend on $n$,
but the ideals $I,I'$ and the threshold $N$ are fixed.
The equality therefore holds for every $n\geq1$,
which is the required equality of Hilbert functions.
\end{proof}

\section{Proofs of the main theorems}\label{sec:main-proof}

We first establish a common estimate under \eqref{eq:koszul-hypothesis},
with $A=R/(\f)$, $s=\dim A>0$, and $D=\hdeg(A)$.
First replace $R$ by its
completion and all ideals, sequences and modules by their base changes,
keeping the same notation. The definitions of $D$ and $k=k(A)$
through completion leave both invariants unchanged. Next pass to
$R'=\widehat{R[X]_{\m R[X]}}$, where $X$ is an indeterminate.
Its residue field is $(R/\m)(X)$, which is infinite.
These are faithfully flat local maps with extended maximal ideal.
A composition series of a finite-length module extends to one over
the new residue field; therefore all finite lengths, maximal-ideal Hilbert
functions and their dimensions and multiplicities are preserved.
Koszul homology commutes with these flat maps.

To check the deficiency modules, after the first completion choose
a complete regular local presentation $P\twoheadrightarrow R$, with
maximal ideal $\mathfrak t$, and put
$P'=\widehat{P[X]_{\mathfrak tP[X]}}$.
A regular system of parameters of $P$ extends to one of $P'$;
thus $P'$ is regular of the same dimension, and $R'\cong R\otimes_PP'$.
For every finite $R$-module $M$, flat base change of a finite free
resolution over $P$ gives
\[
 \Ext_P^{\dim P-j}(M,P)\otimes_PP'
 \cong \Ext_{P'}^{\dim P-j}(M\otimes_PP',P').
\]
Thus the deficiency modules commute with this base change. Since
dimensions, multiplicities and finite lengths are preserved, induction
on $\dim M$ in \eqref{eq:module-hdeg-definition} gives
$\hdeg_{R'}(M\otimes_RR')=\hdeg_R(M)$. In particular, $D$ and $k$
are unchanged. If the original final quotient contains a field,
so do its completion and its faithfully flat extension; thus the
equicharacteristic hypothesis on $A$ is also preserved.
Every perturbation with entries in $\m^N$ has entries
in $(\m R')^N$ after extension, and equality of the required finite
lengths after extension implies equality before extension. We may
therefore assume that $R$ is complete with infinite residue field.

We use the following reduction bound.

\begin{lemma}[Reduction bound for local rings]
\label{lem:ring-reduction-bound}
Let $(A,\mathfrak n)$ be a complete Noetherian local ring with
infinite residue field. Put $s=\dim A>0$, $D=\hdeg(A)$, and
$\rho=\max\{sD-2s+1,0\}$.
For general $y_1,\ldots,y_s\in\mathfrak n$, the sequence
$y_1,\ldots,y_s$ is a filter regular system of parameters and
$Q=(y_1,\ldots,y_s)A$ satisfies
\[
 \mathfrak n^{\rho+1}=Q\mathfrak n^\rho.
\]
\end{lemma}
\begin{proof}
Choose a general minimal reduction $J_0$ of $\mathfrak n$.
Successive applications of Lemma~\ref{lem:generic-power-hdeg}
with exponent one give $\ell_A(A/J_0)\le D$.
If $\ell_A(A/J_0)=1$, then $J_0=\mathfrak n$, so the reduction
number is zero. Otherwise, \cite[Theorem~2.5]{Vasconcelos2003} gives
\[
 r(\mathfrak n)\le s\ell_A(A/J_0)-2s+1\le\rho,
\]
where $r(\mathfrak n)$ denotes the minimum of the reduction numbers
over all minimal reductions of $\mathfrak n$.
To obtain the asserted equality for general $y_1,\ldots,y_s\in\mathfrak n$, put
$\Bbbk=A/\mathfrak n$ and $G=\operatorname{gr}_{\mathfrak n}(A)$,
so $G_j=\mathfrak n^j/\mathfrak n^{j+1}$. For
$\boldsymbol\ell=(\ell_1,\ldots,\ell_s)\in G_1^s$, consider the map
\[
 \Phi_{\boldsymbol\ell}:G_\rho^{\oplus s}\longrightarrow G_{\rho+1},\qquad
 (u_1,\ldots,u_s)\longmapsto\sum_{i=1}^s\ell_i u_i.
\]
Since $r(\mathfrak n)\le\rho$, some minimal reduction
$J_1=(a_1,\ldots,a_s)A$ satisfies
$\mathfrak n^{\rho+1}=J_1\mathfrak n^\rho$.
Reducing this equality modulo $\mathfrak n^{\rho+2}$ shows that
$\Phi_{\boldsymbol a}$ is surjective, where
$\boldsymbol a=(a_1+\mathfrak n^2,\ldots,a_s+\mathfrak n^2)\in G_1^s$.

First fix bases of $G_\rho$ and $G_{\rho+1}$, using the induced
basis of $G_\rho^{\oplus s}$ to represent $\Phi_{\boldsymbol\ell}$
by a matrix. Next choose a basis of $G_1$ to give coordinates to
each $\ell_i$. Since multiplication $G_1\times G_\rho\to G_{\rho+1}$
is $\Bbbk$-bilinear, the matrix entries are linear in these
coordinates. Surjectivity at $\boldsymbol a$ implies that some
minor $\delta$ of size $\dim_{\Bbbk}G_{\rho+1}$ satisfies
$\delta(\boldsymbol a)\ne0$. Thus $\delta$ is a nonzero polynomial.
Set
\[
 U=\{\boldsymbol\ell\in G_1^s:\delta(\boldsymbol\ell)\ne0\}.
\]
Then $U$ is a nonempty Zariski open subset containing
$\boldsymbol a$, and $\Phi_{\boldsymbol\ell}$ is surjective
for every $\boldsymbol\ell\in U$.

For $\boldsymbol\ell\in U$, lift the $\ell_i$ to any $y_i\in\mathfrak n$
and set $Q=(y_1,\ldots,y_s)A$.
Then $G_{\rho+1}=\sum_i\ell_iG_\rho$, or equivalently,
\[
 \mathfrak n^{\rho+1}=Q\mathfrak n^\rho+\mathfrak n^{\rho+2}.
\]
Applying Nakayama's lemma to the finite $A$-module
$\mathfrak n^{\rho+1}/Q\mathfrak n^\rho$ gives
$\mathfrak n^{\rho+1}=Q\mathfrak n^\rho$.
The intersection of $U$ with the nonempty open subset of $G_1^s$
defining filter regular systems of parameters is therefore nonempty.
Any lifts $y_1,\ldots,y_s\in\mathfrak n$ of a tuple in this
intersection satisfy both required conditions, proving the lemma.
\end{proof}

Return to $A=R/(\f)$, and write $\Bbbk=A/\m A$.
For $0\le j<s$, put $M_j=K^j(A)$, $d_j=\dim M_j$,
$h_j=\hdeg(M_j)$, and $w_j=\binom{s-1}{j}$.
Our zero-module conventions give $d_j=-\infty$ and $h_j=0$ if $M_j=0$.
Thus
\begin{equation}\label{eq:delta-sum}
 \Delta=D-e(A)=\sum_{j=0}^{s-1}w_jh_j.
\end{equation}
Set $\rho=\max\{sD-2s+1,0\}$. By Lemma~\ref{lem:ring-reduction-bound}
applied to $(A,\m A)$, general
$y_1,\ldots,y_s\in\m$ have images forming a filter regular
system of parameters of $A$ and give $Q=(y_1,\ldots,y_s)R$ with
\begin{equation}\label{eq:main-rho}
 \m^{\rho+1}A=Q\m^\rho A.
\end{equation}

In each proof below, we impose finitely many additional nonempty open
conditions on prefixes of this tuple. Pulling them back to
$(\m A/(\m A)^2)^s$ and intersecting with the preceding conditions
still gives a nonempty open set.
The chosen classes lift to $R$ because
$\m/\m^2\to\m A/(\m A)^2$ is surjective.

Suppose that there are integers $a_j\ge0$ and, for every
$0\le j<i\le s$, monic polynomials $P_{j,i}(Z)\in A[Z]$ of degree
$a_j$ such that
\begin{equation}\label{eq:common-annihilators}
 P_{j,i}(y_i)M_j=0,\qquad
 [Z^{a_j-l}]P_{j,i}(Z)\in(y_1,\ldots,y_j)^lA
 \quad(1\le l\le a_j).
\end{equation}
Here $[Z^u]P$ denotes the coefficient of $Z^u$ in $P$.
We take $a_j=0$ and $P_{j,i}=1$ when $M_j=0$.
When $M_j\ne0$ and $d_j=0$, take $a_j=h_j=\ell(M_j)$ and
$P_{j,i}(Z)=Z^{a_j}$; indeed, $\m^{a_j}M_j=0$.
Matlis duality is faithful and takes multiplication by an element
to multiplication by the same element. Consequently
$\Ann_A M_j=\Ann_A H^j_{\m}(A)$, so these evaluations also
annihilate the required local cohomology modules.

For $1\le l\le a_j$, lift the coefficient
$[Z^{a_j-l}]P_{j,i}(Z)$ to $(y_1,\ldots,y_j)^lR$,
keeping the leading coefficient equal to $1$.
Retain the notation $P_{j,i}$ for the lifted polynomial in $R[Z]$.
Proposition~\ref{prop:triangular-parameters}
gives elements $z_1,\ldots,z_s\in R$, represented by homogeneous
polynomials in the generators of $Q$ of common degree
\[
 t=\max\left\{1,\sum_{j=0}^{s-1}w_ja_j\right\}.
\]
Their images form a $d$-sequence in $A$, and
$\operatorname{reg}\operatorname{gr}_{\q}(A)=0$ for
$\q=(z_1,\ldots,z_s)R\subseteq Q^t$. Moreover, the identity
\[
 Q^{s(t-1)+1}=\q Q^{(s-1)(t-1)}
\]
holds in $R$. Proposition~\ref{thm:tr-criterion}
therefore gives equality of the $\mathfrak m$-adic Hilbert functions
for all perturbations in $\m^N$, where
$N=\rho+(s+1)(t-1)+3$.
Using the definitions of $\rho$ and $t$, we obtain
the following common bound under \eqref{eq:common-annihilators}:
\begin{equation}\label{eq:common-perturbation-bound}
 p_{\m}(\f)\le
 \max\{sD-2s+1,0\}
 +(s+1)\max\left\{\sum_{j=0}^{s-1}w_ja_j-1,0\right\}+3.
\end{equation}
It remains to construct the polynomials in
\eqref{eq:common-annihilators} and bound their degrees.

\begin{proof}[Proof of Theorem~\ref{thm:equicharacteristic}]
For each nonzero $M_j$ of positive dimension, consider the
nonzero positive-dimensional factors $C$ of its dimension
filtration. Put $h=\dim C$ and $n_C=h!\hdeg(C)$.
By Corollary~\ref{cor:module-reduction-bound}, the condition
\[
 \m^{n_C}C=(y_1,\ldots,y_h)\m^{n_C-1}C
\]
holds on a nonempty open subset $U_C\subseteq V^h$, where
$V=\m A/(\m A)^2$.
There are finitely many such factors, since $0\le j<s$ and each
dimension filtration is finite.
For each factor of dimension $h$, the projection
$\pi_h:V^s\to V^h$ onto the first $h$ coordinates is surjective,
so $\pi_h^{-1}(U_C)$ is a nonempty open subset of $V^s$.
The intersection of these finitely many open subsets with the
nonempty open subsets for \eqref{eq:main-rho} and the filter regular
parameter condition is therefore a nonempty open subset of $V^s$.
Any lifts $y_1,\ldots,y_s\in\m$ of a tuple in this intersection
give all the required reductions simultaneously.
Applying Proposition~\ref{prop:equicharacteristic-annihilator} over $A$
to each such $M_j$, and using
Lemma~\ref{lem:arithmetic-homological-degree}, we obtain monic
polynomials $P_{j,i}$ whose evaluations annihilate $M_j$, with
\begin{equation}\label{eq:a-j-equicharacteristic}
 a_j=\adeg(M_j)\le h_j\qquad(0\le j<s).
\end{equation}
This includes the zero and finite-length modules with the choices
already made. For $d_j>0$, the coefficient $[Z^{a_j-l}]P_{j,i}(Z)$ belongs to
$(y_1,\ldots,y_{d_j})^lA\subseteq(y_1,\ldots,y_j)^lA$
for $1\le l\le a_j$, since $d_j\le j$. Thus \eqref{eq:common-annihilators} holds, and
\eqref{eq:delta-sum} gives
\[
 \sum_{j=0}^{s-1}w_ja_j\le\sum_{j=0}^{s-1}w_jh_j=\Delta.
\]
The common estimate \eqref{eq:common-perturbation-bound} yields
the refined bound
\begin{equation}\label{eq:equicharacteristic-refined}
 p_{\m}(\f)\le
 \max\{sD-2s+1,0\}+(s+1)\max\{D-e(A)-1,0\}+3.
\end{equation}
Since $s,D\ge1$ and $0\le\Delta\le D$, we have
$\rho\le sD-1$ and $\max\{\Delta-1,0\}\le D-1$.
Thus \eqref{eq:equicharacteristic-refined} gives
\[
 p_{\m}(\f)\le sD-1+(s+1)(D-1)+3\le(2s+1)D.
\]
This proves \eqref{eq:equicharacteristic-main}.
\end{proof}

\begin{proof}[Proof of Theorem~\ref{thm:main}]
For each $M_j$ of positive dimension, put
\[
 U_j=M_j/H^0_{\m}(M_j),\qquad n_j=d_j!\hdeg(U_j).
\]
By Corollary~\ref{cor:module-reduction-bound}, each condition
\[
 \m^{n_j}U_j=(y_1,\ldots,y_{d_j})\m^{n_j-1}U_j
\]
holds on a nonempty open subset of $(\m A/(\m A)^2)^{d_j}$.
As in the proof of Theorem~\ref{thm:equicharacteristic}, choose
$y_1,\ldots,y_s\in\m$ so that these finitely many conditions hold
simultaneously with \eqref{eq:main-rho} and the filter regular
parameter condition.
Applying Proposition~\ref{prop:monic-module-annihilator} over $A$
to each such $M_j$, we obtain monic polynomials $P_{j,i}$ whose
evaluations annihilate $M_j$ and which satisfy
\[
 a_j\le(d_j!)^{d_j-1}h_j^{d_j},\qquad
 [Z^{a_j-l}]P_{j,i}(Z)\in(y_1,\ldots,y_{d_j})^lA
 \quad(1\le l\le a_j).
\]
Since $d_j\le j$, the coefficient condition in
\eqref{eq:common-annihilators} follows. Also $d_j\le k$ and
$h_j\ge1$ for nonzero $M_j$, so the monotonicity of
$(d!)^{d-1}$ for $d\ge1$ gives
\begin{equation}\label{eq:a-j-bound}
 a_j\le c_kh_j^k,\qquad c_k=(k!)^{k-1}.
\end{equation}
This inequality also holds for the zero and finite-length modules
with the choices already made.
Because the $w_j$ are positive integers and $k\ge1$, we obtain
\begin{equation}\label{eq:t-estimate}
 \sum_jw_ja_j\le c_k\sum_jw_jh_j^k
 \le c_k\left(\sum_jw_jh_j\right)^k
 =c_k\Delta^k.
\end{equation}
Here all sums over $j$ run from $0$ to $s-1$.
Indeed, $w_jh_j^k\le(w_jh_j)^k$ since $w_j\ge1$, and the
expansion of $(\sum_jw_jh_j)^k$ contains each $(w_jh_j)^k$
together with nonnegative mixed terms.
The common estimate \eqref{eq:common-perturbation-bound} gives
\begin{equation}\label{eq:refined-main}
 p_{\m}(\f)\le
 \max\{sD-2s+1,0\}+(s+1)\max\{c_k\Delta^k-1,0\}+3.
\end{equation}
Since $s,D,k,c_k\ge1$ and $0\le\Delta\le D$, we have
$\rho\le sD-1$, $D\le D^k$, and
$\max\{c_k\Delta^k-1,0\}\le c_kD^k-1$.
Thus \eqref{eq:refined-main} gives
\[
 p_{\m}(\f)\le sD-1+(s+1)(c_kD^k-1)+3\le C_{s,k}D^k.
\]
This proves the theorem.
\end{proof}

\section{Consequences and a comparison example}\label{sec:consequences}

We retain the notation of Theorem~\ref{thm:main}. The refined bound
\eqref{eq:refined-main}, with $\Delta=D-e(A)$ and $c_k=(k!)^{k-1}$,
yields the first two corollaries below without a restriction on
characteristic. We then obtain a linear bound for certain
mixed-characteristic quotients from
Corollary~\ref{cor:prime-power-annihilator}, and conclude with a
comparison of both main bounds with the homological-degree bound
of Quy and Ngo Viet Trung.

\begin{corollary}\label{cor:linear}
Under the hypotheses of Theorem~\ref{thm:main}, if $\dim K^i(\widehat A)\le1$ for every $0\le i<s$, then
\[
 p_{\m}(\f)\le(2s+1)D.
\]
In particular, the linear bound holds whenever $\dim A\le2$.
\end{corollary}
\begin{proof}
The hypothesis gives $k(A)=1$ and $c_1=1$. Thus
\eqref{eq:refined-main} yields the refined bound
\begin{equation}\label{eq:linear-refined}
 p_{\m}(\f)\le
 \max\{sD-2s+1,0\}+(s+1)\max\{D-e(A)-1,0\}+3.
\end{equation}
The asserted linear bound follows from the numerical comparison
in the proof of Theorem~\ref{thm:equicharacteristic}.
For $s\le2$, the inequalities $\dim K^i(\widehat A)\le i$ imply
the hypothesis automatically.
\end{proof}

\begin{corollary}[Generalized Cohen--Macaulay final quotients]
\label{cor:gcm-final-quotient}
Under the hypotheses of Theorem~\ref{thm:main}, suppose that $A=R/(\f)$
is generalized Cohen--Macaulay, meaning that $H^i_{\m}(A)$ has
finite length for every $0\le i<s$. Set
\[
 \tau=\mathbb I(A):=\sum_{i=0}^{s-1}\binom{s-1}{i}
                  \ell_A(H^i_{\m}(A)),\qquad D=e(A)+\tau.
\]
Then
\[
 p_{\m}(\f)\le\max\{sD-2s+1,0\}+(s+1)\max\{\tau-1,0\}+3.
\]
No generalized Cohen--Macaulay assumption is imposed on $R$.
If $A$ is Cohen--Macaulay, the bound becomes
$\max\{se(A)-2s+1,0\}+3$. This agrees with the numerical bound
$s!e(A)+2$ in \cite[Proposition~1]{ST97} when $s=1$, and is
strictly smaller when $s\ge2$.
\end{corollary}
\begin{proof}
Completion preserves these finite lengths. Matlis duality gives
$\ell(K^i(\widehat A))=\ell_A(H^i_{\m}(A))$ for $0\le i<s$, so every
nonzero lower deficiency module has dimension zero and homological
degree equal to that length. Thus $k(A)=1$, and the defining formula
\eqref{eq:module-hdeg-definition} gives $\hdeg(A)=e(A)+\tau$.
The bound follows from \eqref{eq:linear-refined}.
If $A$ is Cohen--Macaulay, all the lower local cohomology modules
vanish, so $\tau=0$.
For a generalized Cohen--Macaulay final quotient, the construction
is explicit:
write $h_j=\ell_A(H^j_{\m}(A))=\ell(K^j(\widehat A))$ for $0\le j<s$
and take $P_{j,i}(Z)=Z^{h_j}$ for $0\le j<i\le s$. Proposition~\ref{prop:triangular-parameters}
then gives
\[
 t=\max\{1,\tau\},\qquad z_i=y_i^t\quad(1\le i\le s).
\]
Thus in this case the triangular parameter construction reduces to
taking a common power of the chosen parameters.
\end{proof}

\begin{corollary}[Linear bound in mixed characteristic]
\label{cor:prime-power-deficiencies}
Under the hypotheses of Theorem~\ref{thm:main}, suppose that
$\operatorname{char}(R/\m)=p>0$ and that, for every $0\le j<s$,
some power of $p$ annihilates $K^j(\widehat A)$.
Then
\[
 p_{\m}(\f)\le(2s+1)D.
\]
\end{corollary}

\begin{proof}
Use the completion and residue-field extension from
Section~\ref{sec:main-proof}. The hypotheses on the deficiency
modules are preserved by the flat base change described there.
Retain the notation $M_j$, $h_j$, and $w_j$ of that section.
Corollary~\ref{cor:prime-power-annihilator} provides the required
monic polynomials for the positive-dimensional $M_j$, with
$a_j=\adeg(M_j)\le h_j$. Their coefficient conditions use the
prefix $y_1,\ldots,y_{d_j}$, where $d_j=\dim M_j\le j$.
The finitely many general choices are compatible with the common
construction. For zero and finite-length modules, use the choices
already made in that construction. Thus
\[
 \sum_{j=0}^{s-1}w_ja_j\le\sum_{j=0}^{s-1}w_jh_j=D-e(A).
\]
Substitution into \eqref{eq:common-perturbation-bound} gives
the refined bound
\begin{equation}\label{eq:prime-power-refined}
 p_{\m}(\f)\le
 \max\{sD-2s+1,0\}+(s+1)\max\{D-e(A)-1,0\}+3.
\end{equation}
The asserted linear bound follows from the numerical comparison
in the proof of Theorem~\ref{thm:equicharacteristic}.
\end{proof}

If $\widehat A$ is equidimensional, the annihilation hypothesis
in Corollary~\ref{cor:prime-power-deficiencies} is equivalent to
$\widehat A[1/p]$ being Cohen--Macaulay.
Indeed, Matlis duality and \cite[Lemma~4.2]{MQS} identify the
non-Cohen--Macaulay locus of $\widehat A$ with the union of the
supports of its lower deficiency modules.
Since these modules are finite, this union is contained in
$V(p)$ if and only if each is annihilated by a power of $p$.

\begin{example}
\label{ex:mixed-characteristic-linear}
Let $V$ be a complete discrete valuation ring of characteristic
zero with maximal ideal $pV$ and infinite residue field
$F=V/pV$ of characteristic $p>0$. For $s\ge3$, put
\[
 S=V[[x_1,\ldots,x_{s-1}]],\qquad A=S\ltimes S/(p).
\]
Here the underlying $S$-module is $S\oplus S/(p)$, with
multiplication $(a,u)(b,v)=(ab,av+bu)$. Equivalently, with
$P=S[[T]]$, we have
\[
 A\cong P/(T^2,pT).
\]
The ring $A$ is complete of dimension $s$ and has mixed
characteristic. Since $p,T$ is a regular sequence in $P$, the
following is a free resolution of $A$ over $P$:
\[
 0\longrightarrow P
 \xrightarrow{\binom{-p}{T}}P^2
 \xrightarrow{(T^2\ \ pT)}P
 \longrightarrow A\longrightarrow0.
\]
As $\dim P=s+1$, applying $\operatorname{Hom}_P(-,P)$ gives
\[
 K^{s-1}(A)\cong P/(p,T)\cong S/(p),\qquad
 K^j(A)=0\quad(0\le j<s-1).
\]
In particular, the only nonzero lower deficiency module has
dimension $s-1\ge2$ and is annihilated by $p$.
Thus Corollary~\ref{cor:prime-power-deficiencies} applies, whereas
the equicharacteristic hypothesis of
Theorem~\ref{thm:equicharacteristic} and the dimension condition
of Corollary~\ref{cor:linear} do not hold.

The exact sequence $0\to S/(p)\to A\to S\to0$ gives
$e(A)=e(S)=1$, because $\dim S/(p)=s-1<s$.
The module $S/(p)$ is Cohen--Macaulay of multiplicity one, so
$D=\hdeg(A)=2$ and $k(A)=s-1$.
Whenever $A$ occurs as the final quotient under the hypotheses
of Theorem~\ref{thm:main}, the refined bound
\eqref{eq:prime-power-refined} gives the sufficient perturbation
order $4$, independently of $s$.
\end{example}

\begin{example}[Comparison with the bound of Quy and Ngo Viet Trung]
\label{ex:qt-hdeg-comparison}
Let $F$ be an infinite field and let $L\ge1$ be an integer. Put
\[
 B=F[[x,y,z,u]],\qquad E=B/(x),\qquad
 T_L=B/(x,y,z,u^L),\qquad R_L=B\ltimes(E\oplus T_L),
\]
using the idealization notation of
Example~\ref{ex:mixed-characteristic-linear}.
The ring $R_L$ is a complete Noetherian local ring of dimension four,
with maximal ideal $\m_L=\mathfrak n\oplus(E\oplus T_L)$, where
$\mathfrak n=(x,y,z,u)B$.
Multiplication by $u=(u,0)\in R_L$ is injective on $B$ and $E$,
and has kernel $Fu^{L-1}$ on $T_L$. Thus $u$ is filter regular.
With $S=F[[x,y,z]]$ and $E_0=S/(x)$, quotienting the three summands
by $u$ gives
\[
 A=R_L/uR_L\simeq S\ltimes(E_0\oplus F),\qquad \dim A=3,
\]
independently of $L$.

For a finite module $M$ over a local ring $(C,\mathfrak n_C)$,
\[
 (\mathfrak n_C\oplus M)^n
 =\mathfrak n_C^n\oplus\mathfrak n_C^{n-1}M\qquad(n\ge1).
\]
Since $\dim(E\oplus T_L)<\dim B$ and
$\dim(E_0\oplus F)<\dim S$, the additional summands do not affect
the leading Hilbert coefficients. Hence $e(R_L)=e(A)=1$.

Use the regular presentations $P=S[[v,w]]\twoheadrightarrow A$
and $P'=B[[v,w]]\twoheadrightarrow R_L$, sending $v,w$ to the
generators of the respective square-zero summands, and the exact sequences
\[
 0\longrightarrow E_0\oplus F\longrightarrow A\longrightarrow S\longrightarrow0,
 \qquad
 0\longrightarrow E\oplus T_L\longrightarrow R_L\longrightarrow B\longrightarrow0.
\]
The modules $S,E_0,F$ over $P$ are quotients by regular sequences
of lengths $2,3,5$, respectively; for $B,E,T_L$ over $P'$,
the corresponding regular sequences have lengths $2,3,6$.
Their Koszul resolutions and the long exact Ext sequences give
\[
 K^0(A)\simeq F,\qquad K^1(A)=0,\qquad K^2(A)\simeq E_0,
\]
\[
 K^0(R_L)\simeq T_L,\qquad K^1(R_L)=K^2(R_L)=0,
 \qquad K^3(R_L)\simeq E.
\]
Indeed, $\Ext_P^j(S,P)=0$ and $\Ext_{P'}^j(B,P')=0$ for $j>2$,
so the Ext modules of the middle and left terms are isomorphic
in degrees greater than two. Since $E_0,E,F$ have homological degree one and
$T_L$ has length $L$, we obtain
\[
 D=\hdeg(A)=3,\qquad k(A)=2,\qquad \hdeg(R_L)=L+2.
\]
Thus Corollary~\ref{cor:linear} does not apply, whereas
Theorem~\ref{thm:equicharacteristic} applies because $A$ contains $F$.

Here $s=3$, $D=3$, $e(A)=1$, $k=2$, $\Delta=2$ and $\rho=4$.
The refined linear bound \eqref{eq:equicharacteristic-refined}
gives the sufficient order $N_{\mathrm{lin}}=11$;
Theorem~\ref{thm:main} and its refinement \eqref{eq:refined-main}
give the orders $99$ and $35$, respectively.
For the homological-degree bound of Quy and Ngo Viet Trung
\cite[Theorem~4.5]{QT}, take $J=\m_L$. The parameters are
\[
 d=4,\qquad r=1,\qquad h=\max\{d!,r\}=24,\qquad
 D_*:=\max\{\hdeg(R_L),\hdeg(A)\}=L+2.
\]
The expression stated there becomes
\[
 N_{\mathrm{QT}}=2^{h-1}D_*^{\,3d!-2}=2^{23}(L+2)^{70}.
\]
That theorem preserves the associated graded ring, and hence also
the Hilbert function. For Hilbert-function stability, our sufficient
order remains $11$ while the cited bound grows without bound as
$L\to\infty$, although the final quotient $A$ is fixed.
\end{example}

\end{document}